\documentclass[a4paper,12pt]{article}

\newcommand{\red}[1]{{\colorbox{red}{{#1}}}}

\usepackage{color,makeidx,bm,latexsym,amscd,amsmath,amssymb,enumerate,stmaryrd,amsthm,float,graphicx,psfrag,geometry,comment}

\usepackage[all]{xy}
\usepackage{bm}

\usepackage{tikz,epic,eso-pic}

\tikzset{v/.style={
  circle, draw, inner sep=2pt, minimum size=6pt, fill=white}}

\theoremstyle{plain}
\newtheorem{theorem}{Theorem}[section]
\newtheorem{corollary}[theorem]{Corollary}

\theoremstyle{definition}
\newtheorem{definition}[theorem]{Definition}
\newtheorem{remark}[theorem]{Remark}
\newtheorem{example}[theorem]{Example}
\newtheorem{proposition}[theorem]{Proposition}

\def\qed{\hfill $\Box$}

\newcommand{\Div}{\operatorname{Div}}
\newcommand{\order}{\operatorname{order}}
\newcommand{\Img}{\operatorname{Im}}
\newcommand{\Ker}{\operatorname{Ker}}
\newcommand{\ch}{\operatorname{ch}}
\newcommand{\bch}{\operatorname{bch}}
\newcommand{\uch}{\operatorname{uch}}
\newcommand{\codim}{\operatorname{codim}}
\newcommand{\Rect}{\operatorname{Rect}}
\newcommand{\Int}{\operatorname{Int}}
\newcommand{\wtilde}{\widetilde}
\newcommand{\gl}{\operatorname{gl}}
\newcommand{\OS}{\operatorname{OS}}
\newcommand{\Sing}{\operatorname{Sing}}
\newcommand{\Sep}{\operatorname{Sep}}
\newcommand{\Fix}{\operatorname{Fix}}

\newcommand{\rA}{\mathrm{A}}
\newcommand{\rB}{\mathrm{B}}
\newcommand{\rC}{\mathrm{C}}
\newcommand{\rD}{\mathrm{D}}
\newcommand{\rE}{\mathrm{E}}
\newcommand{\rF}{\mathrm{F}}
\newcommand{\rG}{\mathrm{G}}
\newcommand{\rH}{\mathrm{H}}
\newcommand{\rI}{\mathrm{I}}
\newcommand{\rJ}{\mathrm{J}}
\newcommand{\rK}{\mathrm{K}}
\newcommand{\rL}{\mathrm{L}}
\newcommand{\rM}{\mathrm{M}}
\newcommand{\rN}{\mathrm{N}}
\newcommand{\rO}{\mathrm{O}}
\newcommand{\rP}{\mathrm{P}}
\newcommand{\rQ}{\mathrm{Q}}
\newcommand{\rR}{\mathrm{R}}
\newcommand{\rS}{\mathrm{S}}
\newcommand{\rT}{\mathrm{T}}
\newcommand{\rU}{\mathrm{U}}
\newcommand{\rV}{\mathrm{V}}
\newcommand{\rW}{\mathrm{W}}
\newcommand{\rX}{\mathrm{X}}
\newcommand{\rY}{\mathrm{Y}}
\newcommand{\rZ}{\mathrm{Z}}

\newcommand{\frA}{\mathfrak{A}}
\newcommand{\frB}{\mathfrak{B}}
\newcommand{\frC}{\mathfrak{C}}
\newcommand{\frD}{\mathfrak{D}}
\newcommand{\frE}{\mathfrak{E}}
\newcommand{\frF}{\mathfrak{F}}
\newcommand{\frG}{\mathfrak{G}}
\newcommand{\frH}{\mathfrak{H}}
\newcommand{\frI}{\mathfrak{I}}
\newcommand{\frJ}{\mathfrak{J}}
\newcommand{\frK}{\mathfrak{K}}
\newcommand{\frL}{\mathfrak{L}}
\newcommand{\frM}{\mathfrak{M}}
\newcommand{\frN}{\mathfrak{N}}
\newcommand{\frO}{\mathfrak{O}}
\newcommand{\frP}{\mathfrak{P}}
\newcommand{\frQ}{\mathfrak{Q}}
\newcommand{\frR}{\mathfrak{R}}
\newcommand{\frS}{\mathfrak{S}}
\newcommand{\frT}{\mathfrak{T}}
\newcommand{\frU}{\mathfrak{U}}
\newcommand{\frV}{\mathfrak{V}}
\newcommand{\frW}{\mathfrak{W}}
\newcommand{\frX}{\mathfrak{X}}
\newcommand{\frY}{\mathfrak{Y}}
\newcommand{\frZ}{\mathfrak{Z}}

\newcommand{\bA}{\mathbb{A}}
\newcommand{\bC}{\mathbb{C}}
\newcommand{\bD}{\mathbb{D}}
\newcommand{\bF}{\mathbb{F}}
\newcommand{\bG}{\mathbb{G}}
\newcommand{\bH}{\mathbb{H}}
\newcommand{\bK}{\mathbb{K}}
\newcommand{\bP}{\mathbb{P}}
\newcommand{\bQ}{\mathbb{Q}}
\newcommand{\bR}{\mathbb{R}}
\newcommand{\bS}{\mathbb{S}}
\newcommand{\bX}{\mathbb{X}}
\newcommand{\bY}{\mathbb{Y}}
\newcommand{\bZ}{\mathbb{Z}}

\newcommand{\cA}{\mathcal{A}}
\newcommand{\cB}{\mathcal{B}}
\newcommand{\cC}{\mathcal{C}}
\newcommand{\cD}{\mathcal{D}}
\newcommand{\cE}{\mathcal{E}}
\newcommand{\cF}{\mathcal{F}}
\newcommand{\cG}{\mathcal{G}}
\newcommand{\cH}{\mathcal{H}}
\newcommand{\cI}{\mathcal{I}}
\newcommand{\cJ}{\mathcal{J}}
\newcommand{\cK}{\mathcal{K}}
\newcommand{\cL}{\mathcal{L}}
\newcommand{\cM}{\mathcal{M}}
\newcommand{\cN}{\mathcal{N}}
\newcommand{\cO}{\mathcal{O}}
\newcommand{\cP}{\mathcal{P}}
\newcommand{\cQ}{\mathcal{Q}}
\newcommand{\cR}{\mathcal{R}}
\newcommand{\cS}{\mathcal{S}}
\newcommand{\cT}{\mathcal{T}}
\newcommand{\cU}{\mathcal{U}}
\newcommand{\cV}{\mathcal{V}}
\newcommand{\cW}{\mathcal{W}}
\newcommand{\cX}{\mathcal{X}}
\newcommand{\cY}{\mathcal{Y}}
\newcommand{\cZ}{\mathcal{Z}}

\newcommand{\sA}{\mathsf{A}}
\newcommand{\sB}{\mathsf{B}}
\newcommand{\sC}{\mathsf{C}}
\newcommand{\sD}{\mathsf{D}}
\newcommand{\sE}{\mathsf{E}}
\newcommand{\sF}{\mathsf{F}}
\newcommand{\sG}{\mathsf{G}}
\newcommand{\sH}{\mathsf{H}}
\newcommand{\sI}{\mathsf{I}}
\newcommand{\sJ}{\mathsf{J}}
\newcommand{\sK}{\mathsf{K}}
\newcommand{\sL}{\mathsf{L}}
\newcommand{\sM}{\mathsf{M}}
\newcommand{\sN}{\mathsf{N}}
\newcommand{\sO}{\mathsf{O}}
\newcommand{\sP}{\mathsf{P}}
\newcommand{\sQ}{\mathsf{Q}}
\newcommand{\sR}{\mathsf{R}}
\newcommand{\sS}{\mathsf{S}}
\newcommand{\sT}{\mathsf{T}}
\newcommand{\sU}{\mathsf{U}}
\newcommand{\sV}{\mathsf{V}}
\newcommand{\sW}{\mathsf{W}}
\newcommand{\sX}{\mathsf{X}}
\newcommand{\sY}{\mathsf{Y}}
\newcommand{\sZ}{\mathsf{Z}}

\newcommand{\bgu}{\bigcup}
\newcommand{\bga}{\bigcap}

\def\qed{\hfill $\Box$}

\title{Divides with cusps and symmetric links}
\author{Sakumi Sugawara \thanks{Department of Mathematics, Graduate School of Science, Hokkaido University, North 10, West 8, Kita-ku, Sapporo 060-0810, JAPAN, E-mail: sugawara.sakumi.f5@elms.hokudai.ac.jp}}
\date{\today}

\begin{document}
\maketitle

\begin{abstract}
A Divide with cusps is the image of a proper generic immersion from finite intervals and circles into a $2$-disk which allows to have cusps. A divide with cusps is the generalization of the notion of the divide which is introduced by A'Campo. From a divide with cusps, we can define the associated link in $S^3$. In this paper, we give the characterization of the link in $S^3$ which can be described as the associated link of a divide with cusps. In particular, we prove that every strongly invertible link and $2$-periodic link can be described as the link of a divide with cusps.
\end{abstract}

\renewcommand{\thefootnote}{\fnsymbol{footnote}}
\renewcommand{\thefootnote}{\arabic{footnote}}

\section{Introduction}
A divide is the image of a proper generic immersion from finite intervals and circles into a 2-disk. The notion of the divide is introduced by A'Campo as the tool to describe the real morsification of plane curve singularity and its generalization \cite{acampo-real,acampo-gen}.
From a divide, we can define the associated link in $S^3$ by corresponding to the unit tangent vector of a divide. 
It is known that if a divide is constructed by the real morsification of a plane curve singularity, then the associated link is equivalent to the algebraic link of the curve singularity. 
There are some constraints for links that are described as the associated link of divides. For example, a divide link has fiberedness \cite{acampo-gen} and quasipositivity \cite{kaw-quasi}.
Thus, there are not many links that can be described as the divide link.
However, we can calculate the monodromy of the fibered link and the Seifert genus of a divide link from the combinatorial information of the divide. Therefore, a divide is a very nice tool to describe links. 
Recently, divides are applied to study $3$-, $4$-dimensional topology \cite{ko-furu, ish-nao}.

Divides have some generalizations, for example, oriented divides \cite{gi-osaka}, free divides \cite{gi-top,ish-nao}, graph divides \cite{kaw-graph}, (ordered Morse) signed divides \cite{cou-per,cou}. 
Divides are also generalized to immersed curves in a compact oriented surface \cite{ish}.
In this paper, we deal with the divide with cusps, which is introduced by the author and Yoshinaga \cite{sug-yos}. A divide with cusps is the generalization of the divide that allows to have cusps. The notion of divide with cusps was introduced in the context of the handle decomposition of the complement of complexified real line arrangements.
They proved that the Kirby diagram of the complement of a complexified real line arrangement is described as the associated link of a certain divide with cusps.
From a divide with cusps, we can define the associated link in $S^3$ similarly to the original divide. 
However, cusps in the divides bring the associated link to having half twists. Thus, the associated link of divides with cusps may be complicated in general. 
In fact, there are non-fibered links that can be described as the associated link of divides with cusps.
Therefore, the class of links which is described as the link of a divide with cusps is strictly larger than the class of links of original divides.
In this paper, we discuss how large the class of links of divides with cusps is.
We succeeded in giving the characterization of the links of divides with cusps.
The following is the main theorem in this paper.

\begin{theorem}
A link $L$ in $S^3$ is described as the associated link of a divide with cusps if and only if there exists an involution $j$ on $S^3$ which has a non-empty fixed point set such that $j(L) = L$. 
\end{theorem}

Here, an \textit{involution} on $S^3$ means an orientation preserving self-diffeomorphism of $S^3$ satisfying $j^2 = id$ and $j \neq id$.
Let us consider when $L=K$ is a knot. Then, $K$ satisfies $j (K) = K$ if and only if $K$ is either a strongly invertible knot or a $2$-periodic knot. We have the following theorem as the special case of the main result.

\begin{theorem}
A knot $K$ is described as the associated link of a divide with cusps if and only if $K$ is either a strongly invertible knot or $2$-periodic knot.
\end{theorem}


Couture proved that every strongly invertible link is described as the link of a signed divide \cite{cou}. A signed divide can be considered as a divide with cusps. Equipping a negative sign to a double point corresponds to adding some cusps around the double point (Remark \ref{rem:negative}). Thus, our main result is the generalization of the Couture's result. 
However, the proof of our main result is based on the idea of Couture's result.
In the proof, we construct a well-projection of a link to the disk and perform an appropriate isotopy to the link of a divide with cusps using this projection.

This paper is organized as follows. In Section \ref{sec:divides}, we introduce the theory of divides with cusps. In Section \ref{sec:involution}, we recall the symmetry and involutions on links. In Section \ref{sec:main}, we prove the main theorem and give some examples, including torus knots and $2$-bridge knots.

\vspace{3mm}
\textbf{Acknowledgement.}
The author would like to thank Professor Masahiko Yoshinaga for the helpful discussions on this research.
The author also would like to thank Professors Mikami Hirasawa, Tomomi Kawamura, and Makoto Sakuma for the discussions during the conference Branched Coverings, Degenerations, and Related Topics (March 2022, online).
This work is supported by the JSPS KAKENHI Grant Number JP22KJ0114 and JP23H00081 (PI: Masahiko Yoshinaga).

\section{Divides with cusps and associated links}
\label{sec:divides}

\begin{figure}[htbp]
\centering

\begin{tikzpicture}

\draw [thick](0,0) circle (2);
\draw (-1,0) .. controls ++(0,0.3) and ++(-0.5, -0.5) .. (0,1);
\draw (1,0) .. controls ++(0,0.3) and ++(0.5, -0.5) .. (0,1);
\draw (0,1) .. controls ++(0.7,0.7) and ++(-0.7,0.7) .. (0,1);
\draw (-1,0) .. controls ++(0,-1) and ++(0,0.7) .. (0,-1);
\draw (1,0) .. controls ++(0,-1) and ++(0,0.7) .. (0,-1);

\draw (1.2,1.6) .. controls ++(-1,-1) and ++(0.7,0) .. (0,0);
\draw (1.2,-1.6) .. controls ++(-1,1) and ++(0.7,0) .. (0,0);


\draw (-1,0.3) circle (0.6);

\draw[densely dotted,thick] (0,-1.4) -- (0,-0.4);
\draw (0,-1) node[left]{$x$};
\fill[black] (0,-1) circle (0.06) ;

\draw (0.2,-1.4)node[below] {$T_x P$};

\end{tikzpicture}

\caption{An example of a divide with cusps}
\label{fig:dividewithcusps}
\end{figure}
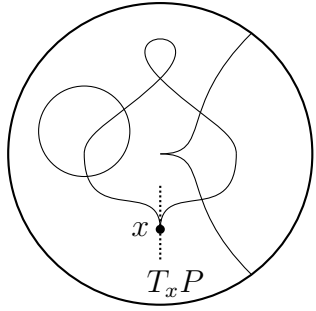

\begin{definition}
Let $X$ be a compact $1$-dimensional manifold.
A \textit{divide with cusps} $P$ is the image of a continuous map $\alpha : (X, \partial X) \rightarrow (D^2, \partial D^2)$ satisfying the following conditions (see Figure \ref{fig:dividewithcusps}).
\begin{itemize}
\item[(i)] $\alpha^{-1} (\partial D^2) = \partial X$, the restriction $\alpha |_{\partial X}$ is injective, and $P$ intersects to $\partial D^2$ transversely.
\item[(ii)] There are finitely many points $p_{1} ,\cdots, p_{s} \in \Int (X)$ such that $\alpha$ is an immersion on $ X \setminus \{p_1 \cdots, p_s\}$ and $\alpha (p_i)$ is a cusp of $P$. A cusp is not a self intersection.
\item [(iii)] There are no singularities except for cusps and double points.
\end{itemize}
\end{definition}

A $1$-dimensional manifold $X$ can be expressed as a disjoint union of finite intervals and circles $\bigsqcup_{j} I_{j} \sqcup \bigsqcup_{k} S_{k}$, where $I_j$ is an interval and $S_k$ is a circle.

\begin{definition}
The image of each connected component of $X$ is called a \textit{component} of a divide with cusps. The image $\alpha(I_{j})$ of each interval is called a \textit{interval component} and the image $\alpha(S_{k})$ of each circle is called a \textit{circle component}.
\end{definition}

\begin{remark}
A divide with cusps that has no cusps is called just a \textit{divide}. 
This is nothing but a divide introduced by A'Campo \cite{acampo-gen}.
\end{remark}

We identify $D^2$ the unit disk in the complex plane $\{z \in \bC \mid |z| \leq 1\}$, and fix the canonical orientation.
Let $TD^2 \cong D^2 \times \bR^2$ be the total space of the tangent bundle of $D^2$.We identify the $3$-sphere $S^3$ to the unit sphere $S(D^2)$ in $TD^2$. Here, 
\[
S(D^2) = \{(x,v) \in TD^2 \mid |x|^2 + |v|^2 = 1\}.
\]
Let $\pi: S(D^2) \rightarrow D^2 \, ; \, \pi(x,v) = x$ be the first projection.
We define an involution $j: S(D^2) \rightarrow S(D^2)$ as $j(x,v) = (x,-v)$.
It is clear that $\Fix (j) = \partial D^2 \times \{0\}$.

\begin{definition}
Let $P$ be a divide with cusps. Then, we define the associated link $L(P)$ in $S(D^2)$ as
\[
L(P) = \{(x,v) \in S(D^2) \mid x \in P, v \in T_{x} P\}.
\]
\end{definition}

\begin{remark}
If $x \in {\rm Int}P$, then there are two vectors $\pm v \in T_{x} P$ which are contained in $L(P)$. If $x \in \partial P$, then only $0 \in T_{x} P$ is contained in $L(P)$. 
The tangent vector is determined by the argument. 
\end{remark}

\begin{figure}[htbp]
\centering

\begin{tikzpicture}
\coordinate (P1) at (0,0);
\coordinate (P2) at (5.2,0);
\coordinate (P3) at (10.4,0);


\draw (P1) --++(4,0);
\draw [->,>=stealth, thick, red] (P1) ++(0,0.03) --++(0.4,0);
\draw [->,>=stealth, thick, red] (P1) ++(1.3,0.03) --++(0.4,0);
\draw [->,>=stealth, thick, red] (P1) ++(2.6,0.03) --++(0.4,0);
\draw [->,>=stealth, thick, red] (P1) ++(4,0.03) --++(0.4,0);

\draw [->,>=stealth, thick, blue] (P1) ++(0,-0.03) --++(-0.4,0);
\draw [->,>=stealth, thick, blue] (P1) ++(1.3,-0.03) --++(-0.4,0);
\draw [->,>=stealth, thick, blue] (P1) ++(2.6,-0.03) --++(-0.4,0);
\draw [->,>=stealth, thick, blue] (P1) ++(4,-0.03) --++(-0.4,0);

\draw [red](P1) ++ (0,-2) ++ (0.45,-0.3) --++(4,0);
\draw [densely dashed, blue](P1) ++ (0,-2) ++ (-0.45,0.3) --++(4,0);

\draw (P1) ++(0,-2) circle [x radius =0.5, y radius =0.8];
\draw (P1) ++(0,-1.2) --++(4,0);
\draw (P1) ++(0,-2.8) --++(4,0);
\draw (P1) ++(4,0)++(0,-2.8) arc [x radius =0.5, y radius =0.8, start angle =-90, end angle =90];
\draw [densely dashed](P1) ++(4,0)++(0,-2.8) arc [x radius =0.5, y radius =0.8, start angle =-90, end angle =-270];


\draw (P2) --++(0.5,0)to[out=0,in=-120] ++(1.5,0.4) to[out=-60,in=180] ++(1.5,-0.4) --++(0.5,0);

\draw [->,>=stealth, thick, red] (P2) ++(0,0.03) --++(0.4,0);
\draw [->,>=stealth, thick, red] (P2) ++(1,0.03) --++(0.4,0);
\draw [->,>=stealth, thick, red] (P2) ++(2,0.4) --++(0,0.4);
\draw [->,>=stealth, thick, red] (P2) ++(3,0.03) --++(-0.4,0);
\draw [->,>=stealth, thick, red] (P2) ++(4,0.03) --++(-0.4,0);

\draw [->,>=stealth, thick, blue] (P2) ++(0,-0.03) --++(-0.4,0);
\draw [->,>=stealth, thick, blue] (P2) ++(1,-0.03) --++(-0.4,0);
\draw [->,>=stealth, thick, blue] (P2) ++(2,0.4) --++(0,-0.4);
\draw [->,>=stealth, thick, blue] (P2) ++(3,-0.03) --++(0.4,0);
\draw [->,>=stealth, thick, blue] (P2) ++(4,-0.03) --++(0.4,0);

\draw [red](P2) ++ (0,-2) ++ (0.45,-0.3) --++(0.5,0) to[out=0,in=180] ++(1,1.1);
\draw [densely dashed, red] (P2) ++ (0,-2) ++ (0.45,-0.3) ++(0.5,0) ++(1,1.1) to[out=0,in=180] ++(0.5,-0.5) --++(1.1,0);
\draw [densely dashed, blue](P2) ++ (0,-2) ++ (-0.45,0.3) --++(1.5,0) to[out=0,in=180] ++(1,-1.1);
\draw [blue](P2) ++ (0,-2) ++ (-0.45,0.3) ++(1.5,0) ++(1,-1.1) to[out=0,in=180] ++(0.5,0.5)--++(1.9,0);

\draw (P2) ++(0,-2) circle [x radius =0.5, y radius =0.8];
\draw (P2) ++(0,-1.2) --++(4,0);
\draw (P2) ++(0,-2.8) --++(4,0);
\draw (P2) ++(4,0)++(0,-2.8) arc [x radius =0.5, y radius =0.8, start angle =-90, end angle =90];
\draw [densely dashed](P2) ++(4,0)++(0,-2.8) arc [x radius =0.5, y radius =0.8, start angle =-90, end angle =-270];

\draw (P3) --++(0.5,0)to[out=0,in=120] ++(1.5,-0.4) to[out=60,in=180] ++(1.5,0.4) --++(0.5,0);

\draw [->,>=stealth, thick, red] (P3) ++(0,0.03) --++(0.4,0);
\draw [->,>=stealth, thick, red] (P3) ++(1,0.03) --++(0.4,0);
\draw [->,>=stealth, thick, red] (P3) ++(2,-0.4) --++(0,-0.4);
\draw [->,>=stealth, thick, red] (P3) ++(3,0.03) --++(-0.4,0);
\draw [->,>=stealth, thick, red] (P3) ++(4,0.03) --++(-0.4,0);

\draw [->,>=stealth, thick, blue] (P3) ++(0,-0.03) --++(-0.4,0);
\draw [->,>=stealth, thick, blue] (P3) ++(1,-0.03) --++(-0.4,0);
\draw [->,>=stealth, thick, blue] (P3) ++(2,-0.4) --++(0,0.4);
\draw [->,>=stealth, thick, blue] (P3) ++(3,-0.03) --++(0.4,0);
\draw [->,>=stealth, thick, blue] (P3) ++(4,-0.03) --++(0.4,0);

\draw [red](P3) ++ (0,-2) ++ (0.45,-0.3) --++(0.5,0) to[out=0,in=180] ++(1,-0.5);
\draw [densely dashed, red] (P3) ++ (0,-2) ++ (0.45,-0.3) ++(0.5,0) ++(1,-0.5) to[out=0,in=180] ++(1,1.1) --++(0.6,0);
\draw [densely dashed, blue](P3) ++ (0,-2) ++ (-0.45,0.3) --++(1.5,0) to[out=0,in=180] ++(1,0.5);
\draw [blue](P3) ++ (0,-2) ++ (-0.45,0.3) ++(1.5,0) ++(1,0.5) to[out=0,in=180] ++(1,-1.1)--++(1.4,0);

\draw (P3) ++(0,-2) circle [x radius =0.5, y radius =0.8];
\draw (P3) ++(0,-1.2) --++(4,0);
\draw (P3) ++(0,-2.8) --++(4,0);
\draw (P3) ++(4,0)++(0,-2.8) arc [x radius =0.5, y radius =0.8, start angle =-90, end angle =90];
\draw [densely dashed](P3) ++(4,0)++(0,-2.8) arc [x radius =0.5, y radius =0.8, start angle =-90, end angle =-270];

\end{tikzpicture}

\caption{Smooth curve, upper cusp, lower cusp, and their links}
\label{fig:halftwist}
\end{figure}
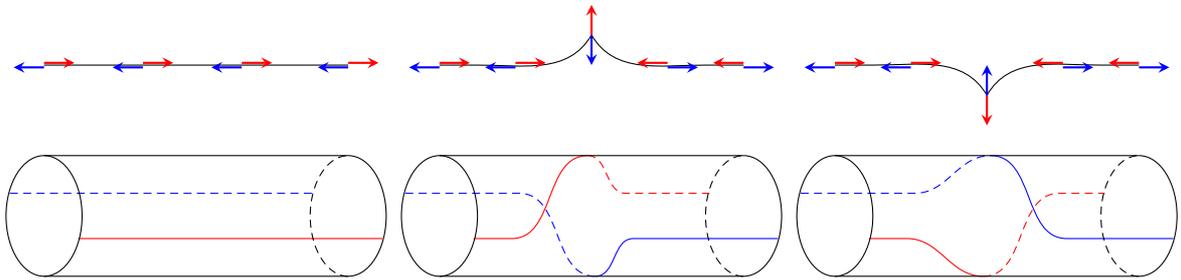

\begin{remark}\label{rem:halftwists}
There is a natural tangent space $T_{x} P$ at a cusp. 
Cusps are corresponding to half twists in the link $L(P)$ locally (see Figure \ref{fig:halftwist}).
\end{remark}

We now introduce a method to describe a link diagram of the associated link of a divide with cusps. The method is based on the algorithm due to Hirasawa \cite{hir}. Let $V \cong D^2 \times S^1$ be a solid torus and $(x,\theta)$ be the coordinate ($x \in D^2$ and $-\pi \leq \theta \leq \pi$).
Let $\tilde{V}$ be the quotient space $V / \sim$ by the equivalence relation $(x_{1}, \theta) \sim (x_{2}, \theta)$ (for $ x_1,x_2 \in \partial D^2$) and $(x,-\pi) \sim (x,\pi)$ (for $x \in D^2$). 
Remark that $\tilde{V}$ is homeomorphic to $S^3$. The map $\phi: \tilde{V} \rightarrow S^3$ defined by $\phi (x, \theta) = (x, \sqrt{1-|u|^2}\cdot e^{\sqrt{-1} \theta})$ gives a homeomorphism. Here, we consider $S^3$ as the unit sphere $S(D^2)$ in the tangent bundle $TD^2$, and tangent space $\bR^2$ as $\bC$. The associated link can be drawn in $\tilde{V}$ by pulling back by this homeomorphism.

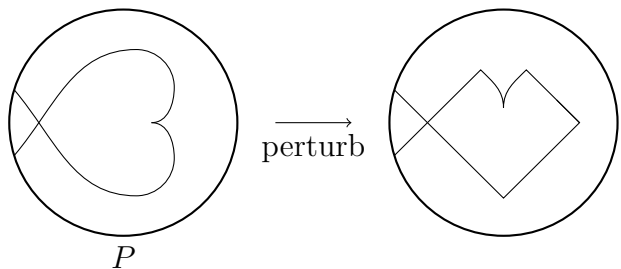
\begin{figure}[htbp]
\centering
\begin{tikzpicture}
\coordinate (P) at (-5,0);
\coordinate (Q) at (0,0);

\draw [thick](P) circle (1.5);
\draw (P) ++(-1.43,0.43) to[out=-50,in=180] ++(1.6,-1.4)to[out=0,in=-90] ++(0.5,0.5) to[out=90,in=0] ++(-0.3,0.47);
\draw (P) ++(-1.43,-0.43) to[out=50,in=180] ++(1.6,1.4)to[out=0,in=90] ++(0.5,-0.5) to[out=-90,in=0] ++(-0.3,-0.47);
\draw (P)++(0,-1.5) node[below]{$P$}; 

\draw [->] (P)++(2,0) --++(1,0);
\draw (P)++(2.5,0) node[below] {perturb};

\draw[thick](Q) circle (1.5);
\draw (Q) ++ (-1,0) --++(0.7,0.7) to[out=-45,in=90] ++(0.3,-0.5) to[out=90,in=-135] ++(0.3,0.5) --++(0.7,-0.7);
\draw (Q) ++ (-1,0) --++ (1,-1) --++(1,1) --++ (-0.3,0.3);
\draw (Q) ++ (-1,0) --++(-0.43,0.43);
\draw (Q) ++ (-1,0) --++(-0.43,-0.43);

\end{tikzpicture}
\caption{Perturbation of a divide with cusps}
\label{fig:perturb}
\end{figure}

Let $P$ be a divide with cusps. By small perturbation, except for the neighborhood of cusps, we assume that $P$ consists of segments with slope $\pm 1$ except for the corners, which is considered as the smooth point where the angle turns $90$ degrees quickly. Around a cusp, we perturb $P$ so that the slope of the cusp is $\pm \infty$ (Figure \ref{fig:perturb}). 
By perturbing $P$ in this way, we can draw the link $L'=\phi^{-1} (L(P))$ in $\tilde{V}$. The link $L'$ consists of linear segments.

\begin{figure}[htbp]
\centering
\begin{tikzpicture}
\coordinate (P) at (0,-2);
\coordinate (P1) at (0,0);
\coordinate (P2) at (0,2);
\coordinate (P3) at (0,4);
\coordinate (P4) at (0,6);

\coordinate (Q) at (-6,3);

\coordinate (R) at (6,3);

\draw (P)circle [x radius=2,y radius=0.6];
\draw (P) ++(-1.63,0.35) --++(2.13,-0.8)--++(1,0.3)--++ (-0.2,0.1);
\draw (P) ++ (-1.5,-0.4) --++ (1.3,0.6) to[out=-0,in=90] ++(0.4,-0.3)to[out=90,in=-180]++(0.5,0.3) --++(0.6,-0.25);

\draw [very thick] (P1) ++ (-1.5,-0.4) --++(1.3,0.6);

\draw[very thick] (P2) ++ (-1.5,-0.4) ++ (1.3,0.6)  ++(0.4,-0.3)++(0.5,0.3) --++(0,2);
\draw[very thick] (P1) ++ (-1.5,-0.4) ++ (1.3,0.6)  ++(0.4,-0.3)++(0.5,0.3) --++(0,-0.5);

\draw[very thick](P2) ++ (-1.5,-0.4) ++ (1.3,0.6) --++(0.3,-0.1) --++(0.3,-1.95) --++(0.3,0.05) ;
\draw[very thick](P4) ++ (-1.5,-0.4) ++ (1.3,0.6) --++(0.3,-0.1)--++(0.3,-1.95) --++(0.3,0.05) ;

\draw [preaction ={draw=white, line width =3 pt},very thick] (P3) ++ (-1.5,-0.4) --++(1.3,0.6);

\draw[preaction={draw = white, line width =3 pt},very thick] (P2) ++ (-1.63,0.35) --++(2.13,-0.8);
\draw[preaction={draw = white, line width =3 pt},very thick] (P4) ++ (-1.63,0.35) --++(2.13,-0.8);

\draw[preaction={draw = white, line width =3 pt},very thick] (P1) ++ (-1.5,-0.4) --++(0,4);

\draw[very thick] (P2) ++ (-1.63,0.35) --++(0,4);


\draw[preaction={draw = white, line width =3 pt},very thick] (P2) ++ (-1.63,0.35) ++(2.13,-0.8) --++(0,2);
\draw[very thick] (P4) ++ (-1.63,0.35) ++(2.13,-0.8) --++(0,0.5);
\draw[very thick] (P1) ++ (-1.63,0.35) ++(2.13,-0.8) --++(0,-0.5);

\draw[very thick](P2) ++ (-1.5,-0.4) ++ (1.3,0.6) --++ (0,2);
\draw[very thick](P4) ++ (-1.5,-0.4) ++ (1.3,0.6) --++ (0,0.5);
\draw[very thick](P1) ++ (-1.5,-0.4) ++ (1.3,0.6) --++ (0,-0.5);

\draw[preaction={draw = white, line width =3 pt},very thick] (P1) ++(-1.63,0.35) ++(2.13,-0.8)--++(1,0.3);
\draw[preaction={draw = white, line width =3 pt},very thick] (P3) ++(-1.63,0.35) ++(2.13,-0.8)--++(1,0.3);

\draw[very thick] (P1) ++(-1.63,0.35) ++(2.13,-0.8) ++(1,0.3) --++ (0,2);
\draw[very thick] (P3) ++(-1.63,0.35) ++(2.13,-0.8) ++(1,0.3) --++ (0,2);

\draw[very thick] (P2) ++ (-1.5,-0.4) ++ (1.3,0.6)  ++(0.4,-0.3)++(0.5,0.3) --++(0.8,-0.35);
\draw[very thick] (P4) ++ (-1.5,-0.4) ++ (1.3,0.6)  ++(0.4,-0.3)++(0.5,0.3) --++(0.8,-0.35);
\draw[very thick] (P4) ++ (-1.5,-0.4) ++ (1.3,0.6)  ++(0.4,-0.3)++(0.5,0.3) --++(0,0.5);

\draw (P1)circle [x radius=2,y radius=0.6];
\draw (P1) ++(-1.63,0.35) --++(2.13,-0.8)--++(1,0.3)--++ (-0.2,0.1);
\draw (P1) ++ (-1.5,-0.4) --++ (1.3,0.6) to[out=-0,in=90] ++(0.4,-0.3)to[out=90,in=-180]++(0.5,0.3) --++(0.6,-0.25);
\draw (P1) ++ (-2,0) node[left]{$-\frac{3}{4} \pi$};

\draw (P2)circle [x radius=2,y radius=0.6];
\draw (P2) ++(-1.63,0.35) --++(2.13,-0.8)--++(1,0.3)--++ (-0.2,0.1);
\draw (P2) ++ (-1.5,-0.4) --++ (1.3,0.6) to[out=-0,in=90] ++(0.4,-0.3)to[out=90,in=-180]++(0.5,0.3) --++(0.6,-0.25);
\draw (P2) ++ (-2,0) node[left]{$-\frac{1}{4} \pi$};

\draw (P3)circle [x radius=2,y radius=0.6];
\draw (P3) ++(-1.63,0.35) --++(2.13,-0.8)--++(1,0.3)--++ (-0.2,0.1);
\draw (P3) ++ (-1.5,-0.4) --++ (1.3,0.6) to[out=-0,in=90] ++(0.4,-0.3)to[out=90,in=-180]++(0.5,0.3) --++(0.6,-0.25);
\draw (P3) ++ (-2,0) node[left]{$\frac{1}{4} \pi$};

\draw (P4)circle [x radius=2,y radius=0.6];
\draw (P4) ++(-1.63,0.35) --++(2.13,-0.8)--++(1,0.3)--++ (-0.2,0.1);
\draw (P4) ++ (-1.5,-0.4) --++ (1.3,0.6) to[out=-0,in=90] ++(0.4,-0.3)to[out=90,in=-180]++(0.5,0.3) --++(0.6,-0.25);
\draw (P4) ++ (-2,0) node[left]{$\frac{3}{4} \pi$};

\draw (P1)++(1.5,1)to[out=0,in=180]++(0.5,-0.5) node[right]{$\phi^{-1} (L(P))$};

\draw (P4) ++ (-2,0.8)node[left]{$\pi$};
\draw (P1) ++ (-2,-0.8)node[left]{$-\pi$};

\draw [preaction={draw = white, line width =3 pt}](P1) ++ (2,0)arc [x radius=2,y radius=0.6, start angle=0, end angle=-90];
\draw [preaction={draw = white, line width =3 pt}](P2) ++ (2,0)arc [x radius=2,y radius=0.6, start angle=0, end angle=-90];
\draw [preaction={draw = white, line width =3 pt}](P3) ++ (2,0)arc [x radius=2,y radius=0.6, start angle=0, end angle=-90];
\draw [preaction={draw = white, line width =3 pt}](P4) ++ (2,0)arc [x radius=2,y radius=0.6, start angle=0, end angle=-90];

\draw[thick] (Q) circle (1.5);
\draw (Q) ++ (-1,0) --++(0.7,0.7) to[out=-45,in=90] ++(0.3,-0.5) to[out=90,in=-135] ++(0.3,0.5) --++(0.7,-0.7);
\draw (Q) ++ (-1,0) --++ (1,-1) --++(1,1) --++ (-0.3,0.3);
\draw (Q) ++ (-1,0) --++(-0.43,0.43);
\draw (Q) ++ (-1,0) --++(-0.43,-0.43);

\draw (Q)++(0,-1.5) node[below]{(perturbed) $P$};


\draw (R) ++(-1.6,-0.5) --++(1.2,1.2);

\draw (R) ++(-1.6,0.5) to[out=135,in=-135] ++ (0,0.2)to[out=45,in=135] ++(0.2,0);
\draw (R) ++(-1.4,-0.7) to[out=-135,in=-45] ++(-0.2,0)to[out=135,in=-135] ++(0,0.2);

\draw[preaction={draw = white, line width =5 pt}] (R) ++(-1.6,0.5) --++(1.3,-1.3);

\draw[preaction={draw = white, line width =5 pt}] (R) ++(-1.4,-0.7) --++(1,1);

\draw[preaction={draw = white, line width =5 pt}] (R) ++(-1.4,0.7) --++(1.3,-1.3);

\draw(R) ++(-1.6,0.5) ++(1.3,-1.3) --++(0.4,-0.4);
\draw[preaction={draw = white, line width =3 pt}] (R) ++(-1.4,0.7) ++(1.3,-1.3)to[out=-45,in=180] ++(0.1,-0.8);

\draw (R) ++(-1.4,-0.7) ++(1,1) --++(0.2,-0.2) to[out=-45,in=-135] ++(0.4,0) --++ (0.5,0.5) --++(0.6,-0.6);
\draw[preaction={draw = white, line width =3 pt}] (R) ++(-1.6,-0.5) ++(1.2,1.2) --++(0.6,-0.6) to[out=-45,in=-135] ++(0.4,0)--++ (0.1,0.1)--++(0.2,-0.2);

\draw (R) ++(-1.4,0.7) ++(1.3,-1.3) ++(0.1,-0.8)--++(0.2,0) to[out=0,in=-135] ++(0.1,0.8) --++(0.6,0.6);
\draw[preaction={draw = white, line width =3 pt}]  (R) ++(-1.6,0.5) ++(1.3,-1.3) ++(0.4,-0.4) --++(0.4,0.4)--++(0.8,0.8);

\draw (R) ++(0,-1.5) node[below]{$L(P)$};

\end{tikzpicture}
\caption{An example of the diagram of $L(P)$}
\label{fig:hirasawa}
\end{figure}
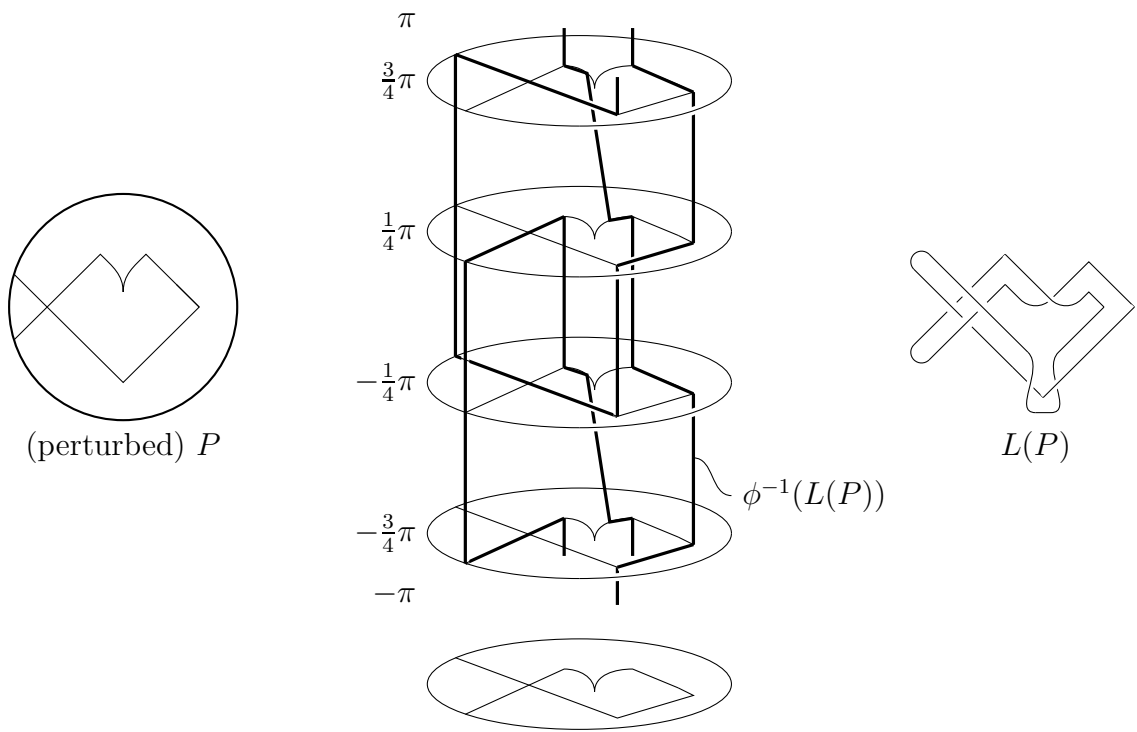

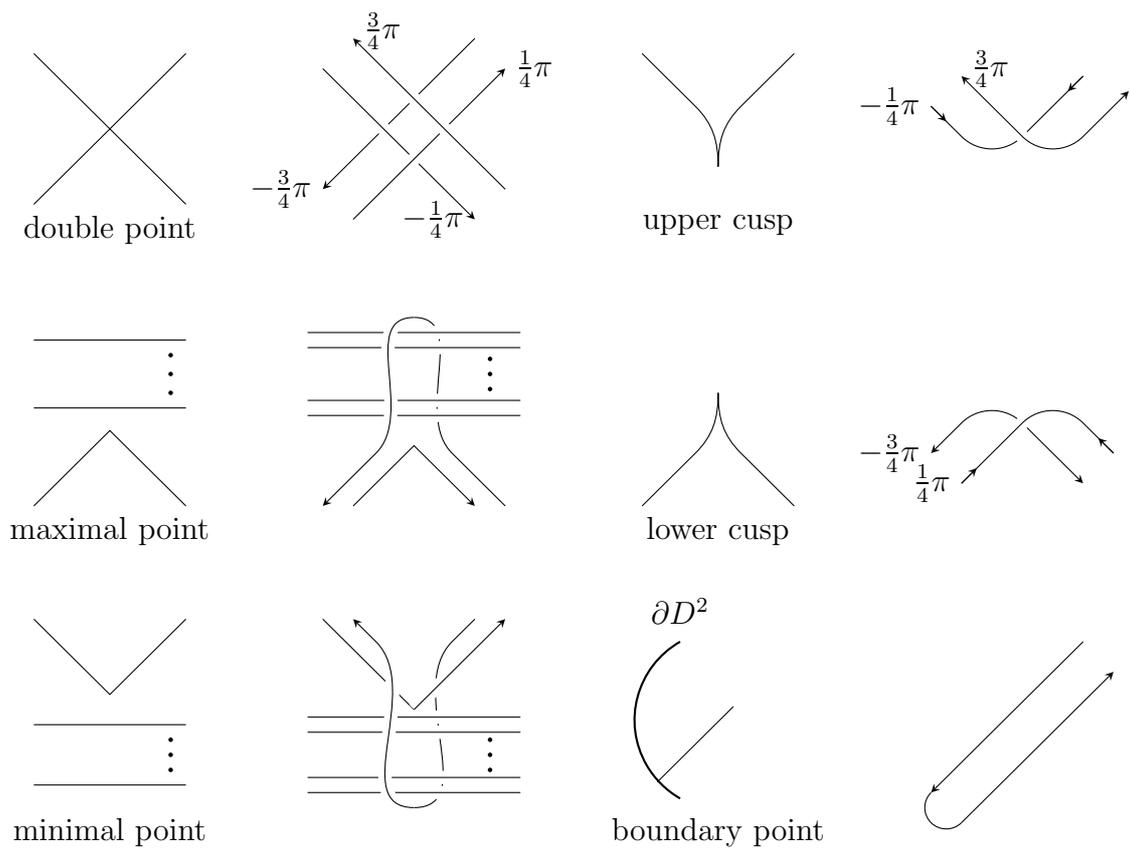
\begin{figure}[htbp]
\centering
\begin{tikzpicture}

\coordinate (P1) at (0,0);
\coordinate (P11) at (4,0);
\coordinate (P2) at (0,-4);
\coordinate (P21) at (4,-4);
\coordinate (P3) at (0,-7.5);
\coordinate (P31) at (4,-7.5);
\coordinate (P4) at (8,0);
\coordinate (P41) at (12,-0.5);
\coordinate (P5) at (8,-4);
\coordinate (P51) at (12,-3.5);
\coordinate (P6) at (8,-8);
\coordinate (P61) at (12,-8); 


\draw(P1) ++ (-1,-1) --++(2,2);
\draw(P1) ++ (-1,1) --++(2,-2);

\draw (P1) ++ (0,-1) node[below]{double point};


\draw[->,>=stealth] (P11) ++ (0.8,1.2) --++(-2,-2) node[left]{$-\frac{3}{4}\pi$};
\draw[draw=white, line width=5 pt] (P11) ++ (-1.2,0.8) --++(2,-2);
\draw[->,>=stealth] (P11) ++ (-1.2,0.8) --++(2,-2)node[left]{$-\frac{1}{4}\pi$};
\draw[draw=white, line width=5pt] (P11) ++ (-0.8,-1.2) --++(2,2); 
\draw[->,>=stealth] (P11) ++ (-0.8,-1.2) --++(2,2)node[right]{$\frac{1}{4}\pi$};
\draw[draw=white, line width=5 pt] (P11) ++ (1.2,-0.8) --++(-2,2) ;
\draw[->,>=stealth] (P11) ++ (1.2,-0.8) --++(-2,2) ;

\draw (P11) ++ (-0.8,1.2)++(0,0.1) node[right]{$\frac{3}{4}\pi$};


\draw (P2) ++(-1,-1) --++(1,1) --++(1,-1);
\draw (P2) ++ (-1,0.3) --++(2,0);
\draw (P2) ++ (-1,1.2) --++(2,0);
\fill (P2) ++(0.8,1) circle (0.03);
\fill (P2) ++(0.8,0.75) circle (0.03);
\fill (P2) ++(0.8,0.5) circle (0.03);

\draw (P2) ++ (0,-1) node[below]{maximal point};

\draw (P21) [->,>=stealth,]++ (-0.8,-1) --++(0.8,0.8)--++(0.8,-0.8);
\draw (P21) ++ (1.2, -1) --++(-0.7,0.7) to[out=135,in=-0] ++(-0.5,1.8);

\draw[preaction={draw = white, line width =5 pt}] (P21) ++ (-1.4,1.1) --++(2.8,0);
\draw [preaction={draw = white, line width =5 pt}](P21) ++ (-1.4,1.3) --++(2.8,0);

\draw[preaction={draw = white, line width =5pt}] (P21) ++ (-1.4,0.2) --++(2.8,0);
\draw[preaction={draw = white, line width =5 pt}] (P21) ++ (-1.4,0.4) --++(2.8,0);

\fill (P21) ++(1,0.95) circle (0.03);
\fill (P21) ++(1,0.75) circle (0.03);
\fill (P21) ++(1,0.55) circle (0.03);

\draw[draw = white, line width =5 pt] (P21) ++ (-1.2,-1) --++(0.7,0.7) to [out=45,in=-180] ++ (0.5,1.8);
\draw[<-,>=stealth] (P21) ++ (-1.2,-1) --++(0.7,0.7) to [out=45,in=-180] ++ (0.5,1.8);


\draw (P3) ++(-1,1) --++(1,-1) --++(1,1);
\draw (P3) ++ (-1,-0.4) --++(2,0);
\draw (P3) ++ (-1,-1.2) --++(2,0);
\fill (P3) ++(0.8,-1) circle (0.03);
\fill (P3) ++(0.8,-0.8) circle (0.03);
\fill (P3) ++(0.8,-0.6) circle (0.03);

\draw (P3) ++ (0,-1.5) node[below]{minimal point};


\draw (P31) ++ (0.8, 1) --++(-0.3,-0.3) to[out=-135,in=-0] ++(-0.5,-2.2);

\draw (P31) ++ (-1.2,1) --++(1.2,-1.2);
\draw [draw = white, line width =5 pt](P31) ++ (-1.2,1) ++(1.2,-1.2)--++(1.2,1.2);
\draw [->,>=stealth](P31) ++ (-1.2,1) ++(1.2,-1.2)--++(1.2,1.2);

\draw[preaction={draw = white, line width =5 pt}] (P31) ++ (-1.4,-1.1) --++(2.8,0);
\draw [preaction={draw = white, line width =5 pt}](P31) ++ (-1.4,-1.3) --++(2.8,0);

\draw[preaction={draw = white, line width =5pt}] (P31) ++ (-1.4,-0.3) --++(2.8,0);
\draw[preaction={draw = white, line width =5 pt}] (P31) ++ (-1.4,-0.5) --++(2.8,0);

\fill (P31) ++(1,-1) circle (0.03);
\fill (P31) ++(1,-0.8) circle (0.03);
\fill (P31) ++(1,-0.6) circle (0.03);

\draw[draw = white, line width =5 pt] (P31)++ (-0.8,1) --++(0.3,-0.3) to [out=-45,in=-180] ++ (0.5,-2.2);
\draw [<-,>=stealth](P31)++ (-0.8,1) --++(0.3,-0.3) to [out=-45,in=-180] ++ (0.5,-2.2);


\draw (P4) ++ (-1,1) --++(0.7,-0.7)to[out=-45,in=90]++(0.3,-0.8)to[out=90,in=-135] ++(0.3,0.8)--++(0.7,0.7); 


\draw(P41) ++ (-1.2,0.8) --++(0.4,-0.4) to[out=-45,in=-135] ++ (0.8,0) --++(0.8,0.8); 
\draw[draw = white, line width = 5 pt] (P41) ++ (-0.8,1.2) --++(0.8,-0.8)to[out=-45,in=-135] ++(0.8,-0) --++(0.6,0.6);
\draw[<->,>=stealth] (P41) ++ (-0.8,1.2) --++(0.8,-0.8)to[out=-45,in=-135] ++(0.8,-0) --++(0.6,0.6);


\draw (P41) ++ (-1.2,0.8) node[left]{$-\frac{1}{4}\pi$};
\draw (P41) ++ (-0.8,1.2)++(0,0.1) node[right]{$\frac{3}{4}\pi$};

\draw[->,>=stealth] (P41) ++ (-1.2,0.8) --++(0.2,-0.2);
\draw[->,>=stealth] (P41) ++ (0.8,1.2) --++(-0.2,-0.2);

\draw (P4) ++(0,-1) node[below] {upper cusp};


\draw (P5) ++ (-1,-1) --++(0.7,0.7) to[out=45,in=-90] ++ (0.3,0.8)to[out=-90,in=135] ++(0.3,-0.8) --++(0.7,-0.7);


\draw[<->,>=stealth] (P51) ++ (-1.2,-0.8) --++(0.4,0.4) to[out=45,in=135] ++(0.8,0)--++(0.8,-0.8);
\draw[preaction={draw = white, line width = 5 pt}](P51) ++ (-0.8,-1.2) --++(0.8,0.8)to[out=45,in=135] ++(0.8,0) --++(0.4,-0.4);

\draw[->,>=stealth] (P51) ++ (-0.8,-1.2) --++(0.2,0.2);
\draw[->,>=stealth] (P51) ++ (1.2,-0.8) --++(-0.2,0.2);

\draw (P51) ++ (-1.2,-0.8) node[left]{$-\frac{3}{4}\pi$};
\draw (P51) ++ (-0.8,-1.2)++(0,0) node[left]{$\frac{1}{4}\pi$};

\draw (P5) ++(0,-1) node[below] {lower cusp};


\draw[thick] (P6) ++ (-0.5,1.2) arc [radius = 1.2, start angle=120, end angle =240] ;
\draw (P6)++(-0.5,1.3) node[above]{$\partial D^2$};
\draw (P6) ++ (-0.8,-0.66) --++(1,1);

\draw (P6) ++ (0,-1) node[below]{boundary point};

\draw[->,>=stealth] (P61) ++(-0.8,-1.2) --++ (2,2);
\draw[->,>=stealth] (P61) ++(0.8,1.2) --++ (-2,-2);
\draw (P61) ++ (-1.2,-0.8) to[out=-135,in=135] ++(-0,-0.4) to[out=-45,in=-135] ++(0.4,0);

\end{tikzpicture}
\caption{Local link diagrams of divides with cusps}
\label{fig:localdiagram}
\end{figure}

The link diagram can be obtained by using the projection $p:(x,\theta) \mapsto x$. After a small perturbation for $L'$, we can draw the link diagram of $L'$ as in Figure \ref{fig:hirasawa}. 
Note that around maximal points, minimal points, and cusps, the angle of the tangent vector crosses $\pm \pi$. 
The link diagrams are locally given as in Figure \ref{fig:localdiagram}.





\begin{remark}\label{rem:moves}
There are moves of divides with cusps that do not change the isotopy type of the link $L(P)$ (Figure \ref{fig:cuspmoves}). We can check them by drawing diagrams.

\begin{figure}[htbp]
\centering 
\begin{tikzpicture}[scale=0.9]
\coordinate (P1) at (0,-3);
\coordinate (P11) at (5,-3);
\coordinate (P2) at (0,-6);
\coordinate (P22) at (5,-6);
\coordinate (P3) at (9.5,0);
\coordinate (P33) at (14.5,0);
\coordinate (P4) at (9.5,-3);
\coordinate (P44) at (14.5,-3);
\coordinate (P5) at (9.5,-6);
\coordinate (P55) at (14.5,-6);
\coordinate (P6) at (0,0);
\coordinate (P66) at (5,0);

\draw (P1) node{(2)};
\draw (P1) ++(1,-2) to[out=30,in=-90] ++(2,1) to[out=90,in=-90] ++(-0.5,1)to[out=-90,in=90]++(-0.5,-1) to[out=-90,in=150] ++(2,-1);
\draw[<->] (P1) ++ (4.5,-1)--++(1,0) ;
\draw (P11) ++(1,-2) to[out=30,in=-90] ++(1.5,0.8)to[out=-90,in=150] ++(1.5,-0.8);

\draw (P2) node {(3)};
\draw [thick] (P2) ++(1.5,-0.5) arc (150:210:2);
\draw (P2) ++(1.25,-1.5) to[out=0,in=-90] ++(1,0.2) to[out=-90,in=180] ++(1,-0.2) ;
\draw (P2) ++(1.34,-0.8) to[out=20,in=200] ++(0.5,0.3)node[right]{$\partial D^2$} ;

\draw[<->] (P2) ++ (4.5,-1.5)--++(1,0) ;

\draw [thick] (P22) ++(1.5,-0.5) arc (150:210:2);
\draw (P22) ++(1.25,-1.5) to[out=0,in=180] ++(1,-0.1) to[out=0,in=180] ++(1,0.1) ;
\draw (P22) ++(1.34,-0.8) to[out=20,in=200] ++(0.5,0.3)node[right]{$\partial D^2$} ;

\draw (P3) node {(4)};
\draw (P3) ++(1,-1.5) to[out=0,in=-90] ++(1.5,1.2) to[out=-90,in=180] ++(1.5,-1.2);
\draw (P3) ++ (1,-0.7) --++(3,0);
\draw [<->](P3) ++ (4.5,-1) --++(1,0);

\draw (P33) ++(1,-1.5) to[out=0,in=-90] ++(1.5,0.7) to[out=-90,in=180] ++(1.5,-0.7);
\draw (P33) ++ (1,-0.7) --++(0.5,0) to[out=0,in=180] ++(1,0.2)to[out=0,in=180] ++(1,-0.2) --++(0.5,0) ;

\draw (P4) node {(5)};
\draw (P4) ++(1.5,-2) --++(2,2);
\draw (P4) ++(1.5,0) --++(2,-2);
\draw (P4) ++(1.5,-1)  to[out=0,in=180] ++(1,0.5) to[out=0,in=180] ++(1,-0.5) ;

\draw [<->](P4) ++ (4.5,-1) --++(1,0);

\draw (P44) ++(1.5,-2) --++(2,2);
\draw (P44) ++(1.5,0) --++(2,-2);
\draw (P44) ++(1.5,-1)  to[out=0,in=180] ++(1,-0.5) to[out=0,in=180] ++(1,0.5) ;

\draw (P5) node {(6)};
\draw (P5) ++ (1.5,-2) --++(2,2);
\draw (P5) ++(1.5,0) to[out=-45,in=45] ++(0.4,-0.8)to[out=45,in=135]++(0.6,-0.2) to[out=-45,in=-135] ++(0.6,-0.2) to[out=-135,in=135] ++(0.4,-0.8);

\draw [<->](P5) ++ (4.5,-1) --++(1,0);

\draw (P55) ++ (1.5,0) --++(2,-2);
\draw (P55) ++ (1.5,-2) to[out=45,in=-45] ++(0.4,0.8)to[out=-45,in=-135]++(0.6,0.2) to[out=45,in=135] ++(0.6,0.2) to[out=135,in=-135] ++(0.4,0.8);

\draw (P6) node {(1)};
\draw (P6) ++(1.5,-1) --++(0.25,0)to[out=0,in=-90] ++(0.4,0.4)to[out=-90,in=180] ++(0.4,-0.4)to[out=0,in=90] ++(0.4,-0.4)to[out=90,in=180] ++(0.4,0.4)--++(0.25,0);

\draw [<->](P6) ++ (4.5,-1) --++(1,0);

\draw (P66) ++(1.5,-1) --++(2,0);

\end{tikzpicture}

\caption{Moves of divides with cusps}
\label{fig:cuspmoves}
\end{figure}
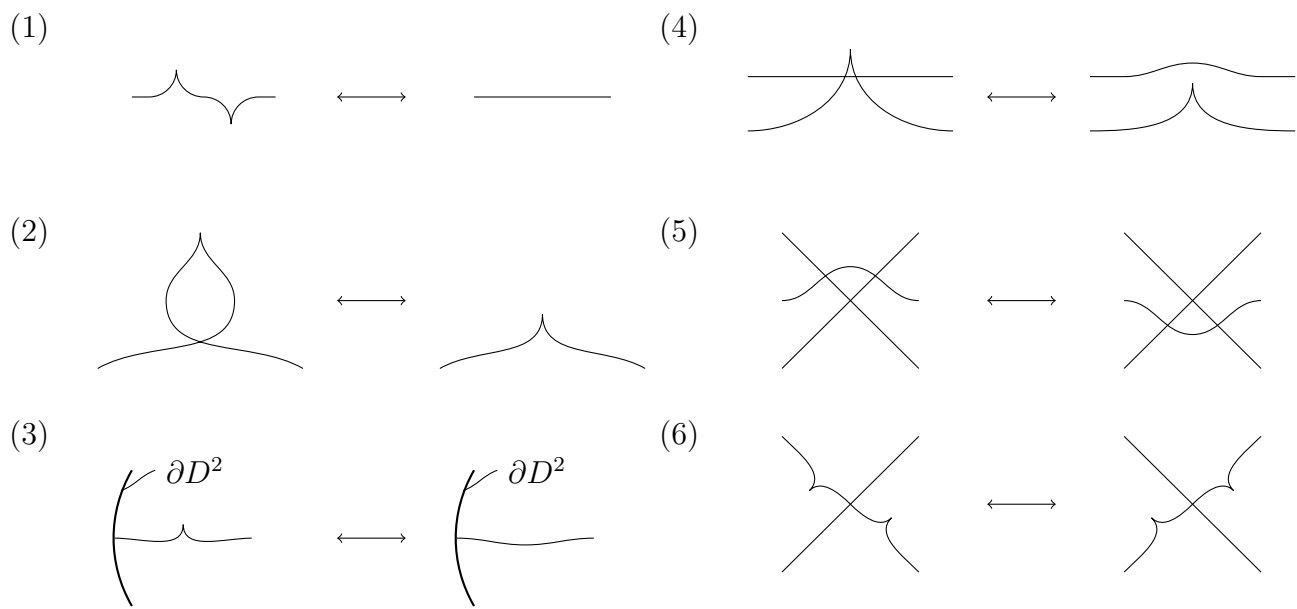

\end{remark}

\section{Involutions on links}
\label{sec:involution}
In this section, we review the symmetry of links and discuss the symmetry of the links of divide with cusps.

\begin{definition}
Let $L$ be an oriented link in $S^3$. A link $L$ is called a \textit{strongly invertible link} if there exists an orientation preserving self-homeomorphism $j:S^3 \rightarrow S^3$ such that 
\begin{itemize}
\item[(i)] $j^2 = id$ and $j \neq id$.
\item[(ii)] The fixed point set $\Fix(j)$ is non-empty set.
\item[(iii)] $j(L) = L$, and the orientation of $j(L)$ is reversed by $j$.
\end{itemize}
We sometimes call the pair $(L,j)$ a strongly invertible link.
\end{definition}

\begin{figure}[htbp]
\centering
\begin{tikzpicture}
\coordinate (P1) at (0,0);
\coordinate (P2) at (7,0);


\draw[thick] (P1) --++(0,-5);
\draw (P1) node[above]{$\Fix (j)$};

\draw[preaction ={draw=white, line width =6 pt}] (P1) ++(0,-4) to[out=-150,in=-30] ++(-2,-1) to [out=150,in=-90] ++(-0.5,1)to[out=90,in=180] ++ (2.5,2);

\draw [preaction ={line width =6 pt, draw=white}](P1) ++ (0,-0.5) to[out=180,in=70] ++(-1.5,-1) to[out=-110,in=150] ++(1.5,-2.5);

\draw[preaction ={draw=white, line width =6 pt}] (P1) ++ (0,-0.5) to[out=0,in=110] ++(1.5,-1) to[out=-70,in=30] ++(-1.5,-2.5);

\draw [preaction ={draw=white, line width =6 pt}](P1) ++(0,-4) to[out=-30,in=-150] ++(2,-1) to [out=30,in=-90] ++(0.5,1)to[out=90,in=0] ++ (-2.5,2);

\draw[preaction ={draw=white, line width =6 pt}] (P1) ++(0.167,-3.9) to[out=-150,in=-30] ++(-2.167,-1.1);

\draw[thick] (P1) --++(0,-3);
\fill[black] (P1) ++(0,-0.5) circle (0.06);
\fill[black] (P1) ++(0,-2) circle (0.06);


\draw[thick] (P2) --++(0,-5);
\draw (P2) node[above]{$\Fix (j)$};

\draw [preaction ={draw=white, line width =6 pt}](P2) ++(0,-1) to[out=-60,in=90] ++(2,-2.5) to [out=-90,in=0] ++(-1,-1.3) to[out=180,in=-60] ++(-1,0.8);

\draw [preaction ={draw=white, line width =6 pt}] (P2) ++ (0,-4) to[out=120,in=-90] ++(-2,2.5) to [out=90,in=180] ++ (1,1.3) to [out=0,in=120] ++ (1,-0.8);

\draw [preaction ={draw=white, line width =6 pt}] (P2) ++ (0,-4) to[out=60,in=-90] ++(2,2.5) to [out=90,in=0] ++ (-1,1.3) to [out=180,in=60] ++ (-1,-0.8);

\draw [preaction ={draw=white, line width =6 pt}](P2) ++(0,-1) to[out=-120,in=90] ++(-2,-2.5) to [out=-90,in=180] ++(1,-1.3) to[out=0,in=-120] ++(1,0.8);

\draw [preaction ={draw=white, line width =6 pt}] (P2) ++ (-0.1,-3.84) to[out=120,in=-90] ++(-1.9,2.36);

\end{tikzpicture}


\caption{Examples of strongly invertible links}
\label{fig:invknot}
\end{figure}
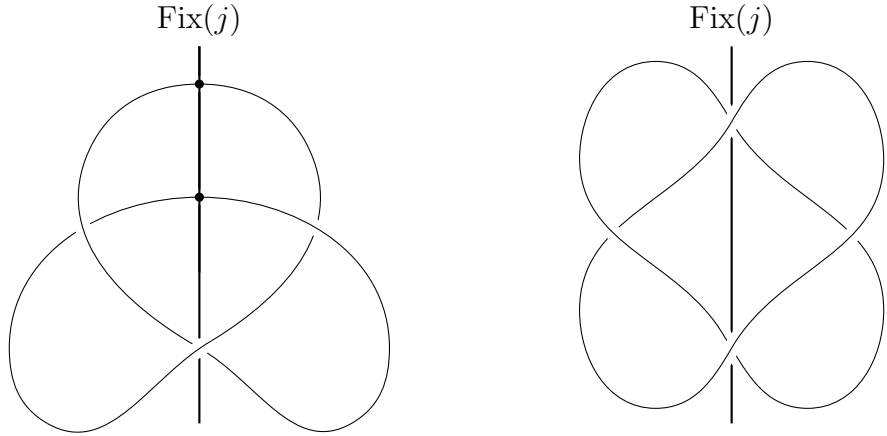

See Figure \ref{fig:invknot} for examples of strongly invertible links. Remark that if $L$ is a strongly invertible knot, then $\Fix (j) \cap L$ is a set consisting of two points.

\begin{definition}
Let $L$ be an oriented link in $S^3$ and $p$ be a prime number. A link $L$ is called \textit{$p$-periodic} if there exists an orientation preserving self-homeomorphism $j:S^3 \rightarrow S^3$ such that
\begin{itemize}
\item[(i)] $j^{p} = id $ and $j^{k} \neq id $ for $1 \leq k \leq p-1$.
\item[(ii)] The fixed point set $\Fix(j)$ is non-empty set.
\item[(iii)] $j(L) = L$, and the orientation of $j(L)$ preserved by $j$.
\end{itemize}
We sometimes call the pair $(L,j)$ is a $p$-periodic link.
\end{definition}

\begin{figure}[htbp]
\centering
\begin{tikzpicture}
\coordinate (P1) at (-7,-2.5);
\coordinate (P2) at (0,0);


\draw (P1) ++ (0,1) ++(0.2,0.2) --++(-0.4,-0.4) to[out=-135, in=0] ++ (-1,-0.2) to[out=180, in= -90]++(-1,2)to[out=90,in=180] ++ (1,2) to[out=0,in=135] ++ (1,-0.4);
\draw (P1) ++ (0,2.5) ++(0.2,0.2) --++(-0.4,-0.4) to[out=-135,in=135]++(-0,-1.1)  ;
\draw (P1) ++ (0,4) ++(0.2,0.2) --++(-0.4,-0.4) to[out=-135,in=135]++(-0,-1.1) ;

\draw[preaction={draw=white, line width =6 pt}, thick] (P1) --++(0,5)node[above]{$\Fix (j)$};

\draw[preaction={draw=white, line width =6 pt}] (P1) ++ (0,1) ++(-0.2,0.2) --++(0.4,-0.4) to[out=-45, in=180] ++ (1,-0.2) to[out=0, in= -90]++(1,2)to[out=90,in=0] ++ (-1,2) to[out=180,in=45] ++ (-1,-0.4);

\draw [preaction={draw=white, line width =6 pt}] (P1) ++ (0,2.5) ++(-0.2,0.2) --++(0.4,-0.4)to[out=-45,in=45]++(-0,-1.1) ;
\draw [preaction={draw=white, line width =6 pt}] (P1) ++ (0,4) ++(-0.2,0.2) --++(0.4,-0.4)to[out=-45,in=45]++(-0,-1.1) ;

\fill (P2) circle (0.08);

\draw (P2)to[out=-90,in=150] ++(0,-0.4) node[right]{$\Fix (j)$};

\draw[->] (P2) ++ (1.9,1.9) to[out=100,in=0] ++(-0.5,0.5);
\draw (P2) ++ (1.9,2.5) node{$j$};



\draw (P2) ++ (0,1.6) to[out=180,in=40] ++ (-0.3,-0.1) to[out=220,in=40] ++ (-0.9,0.2)to[out=220,in=60] ++(-0,-0.8)to[out=240,in=72] ++(-0.3,-0.4); 

\draw[preaction={draw = white, line width =6 pt}] (P2) ++ (0,2.3) to[out=180,in=40] ++ (-0.7,-0.3) to[out=220,in=40] ++(-0.1,-0.8) to [out=220,in=30] ++(-0.9,0.1) to[out=210,in=72] ++(-0.5,-0.6) ;

\draw[preaction={draw = white, line width =6 pt}] (P2) ++ (0,1.6) ++ (-0.3,-0.1)++ (-0.9,0.2)to[out=220,in=60] ++(-0,-0.8)to[out=240,in=72] ++(-0.3,-0.4); 

\draw[rotate =72] (P2) ++ (0,1.6) to[out=180,in=40] ++ (-0.3,-0.1) to[out=220,in=40] ++ (-0.9,0.2)to[out=220,in=60] ++(-0,-0.8)to[out=240,in=72] ++(-0.3,-0.4); 

\draw[rotate = 72, preaction={draw = white, line width =6 pt}] (P2) ++ (0,2.3) to[out=180,in=40] ++ (-0.7,-0.3) to[out=220,in=40] ++(-0.1,-0.8) to [out=220,in=30] ++(-0.9,0.1) to[out=210,in=72] ++(-0.5,-0.6) ;

\draw[rotate=72, preaction={draw = white, line width =6 pt}] (P2) ++ (0,1.6) ++ (-0.3,-0.1)++ (-0.9,0.2)to[out=220,in=60] ++(-0,-0.8)to[out=240,in=72] ++(-0.3,-0.4);

\draw[rotate =144] (P2) ++ (0,1.6) to[out=180,in=40] ++ (-0.3,-0.1) to[out=220,in=40] ++ (-0.9,0.2)to[out=220,in=60] ++(-0,-0.8)to[out=240,in=72] ++(-0.3,-0.4); 

\draw[rotate = 144, preaction={draw = white, line width =6 pt}] (P2) ++ (0,2.3) to[out=180,in=40] ++ (-0.7,-0.3) to[out=220,in=40] ++(-0.1,-0.8) to [out=220,in=30] ++(-0.9,0.1) to[out=210,in=72] ++(-0.5,-0.6) ;

\draw[rotate=144, preaction={draw = white, line width =6 pt}] (P2) ++ (0,1.6) ++ (-0.3,-0.1)++ (-0.9,0.2)to[out=220,in=60] ++(-0,-0.8)to[out=240,in=72] ++(-0.3,-0.4); 

\draw[rotate =216] (P2) ++ (0,1.6) to[out=180,in=40] ++ (-0.3,-0.1) to[out=220,in=40] ++ (-0.9,0.2)to[out=220,in=60] ++(-0,-0.8)to[out=240,in=72] ++(-0.3,-0.4); 

\draw[rotate = 216, preaction={draw = white, line width =6 pt}] (P2) ++ (0,2.3) to[out=180,in=40] ++ (-0.7,-0.3) to[out=220,in=40] ++(-0.1,-0.8) to [out=220,in=30] ++(-0.9,0.1) to[out=210,in=72] ++(-0.5,-0.6) ;

\draw[rotate=216, preaction={draw = white, line width =6 pt}] (P2) ++ (0,1.6) ++ (-0.3,-0.1)++ (-0.9,0.2)to[out=220,in=60] ++(-0,-0.8)to[out=240,in=72] ++(-0.3,-0.4);

\draw[rotate =288] (P2) ++ (0,1.6) to[out=180,in=40] ++ (-0.3,-0.1) to[out=220,in=40] ++ (-0.9,0.2)to[out=220,in=60] ++(-0,-0.8)to[out=240,in=72] ++(-0.3,-0.4); 

\draw[rotate = 288, preaction={draw = white, line width =6 pt}] (P2) ++ (0,2.3) to[out=180,in=40] ++ (-0.7,-0.3) to[out=220,in=40] ++(-0.1,-0.8) to [out=220,in=30] ++(-0.9,0.1) to[out=210,in=72] ++(-0.5,-0.6) ;

\draw[rotate=288, preaction={draw = white, line width =6 pt}] (P2) ++ (0,1.6) ++ (-0.3,-0.1)++ (-0.9,0.2)to[out=220,in=60] ++(-0,-0.8)to[out=240,in=72] ++(-0.3,-0.4);

\draw[densely dashed] (P2) --++(0,2.5);
\draw[densely dashed, rotate =72] (P2) --++(0,2.5);
\draw[densely dashed, rotate =144] (P2) --++(0,2.5);
\draw[densely dashed, rotate =216] (P2) --++(0,2.5);
\draw[densely dashed, rotate =288] (P2) --++(0,2.5);
\end{tikzpicture}
\caption{Examples of periodic links}
\label{fig:periodic}
\end{figure}
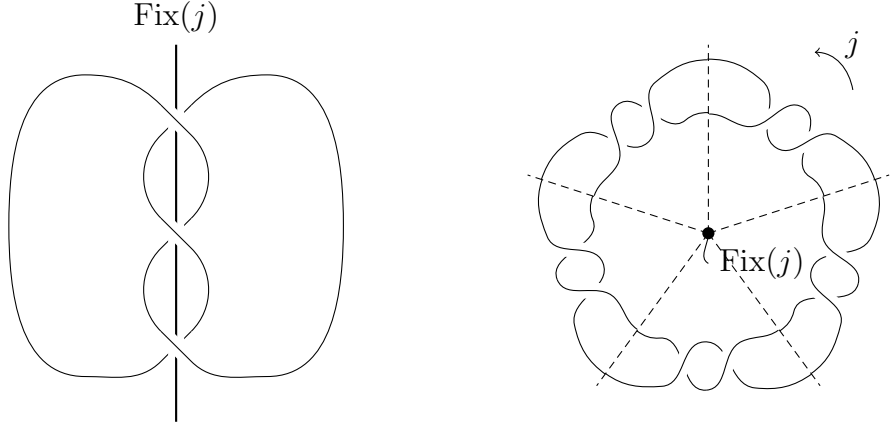

See Figure \ref{fig:periodic} for the examples of $p$-periodic links.
In this paper, we mainly deal with $2$-periodic links.

\begin{remark}
By the resolved Smith's conjecture, it is known that if $j$ is the orientation preserving self-homeomorphism of finite orders with non-empty fixed point sets, then the fixed point set is a trivial knot in $S^3$ \cite{mor-bas}.
\end{remark}

\begin{definition}
Let $(L_{i}, j_{i})$ be a pair of a link and an involution on $S^3$ such that $j_{i} (L_{i}) = L_{i}$ ($i=1,2$).
We call $(L_1,j_1)$ and $(L_2, j_2)$ are \textit{strongly equivalent} if there there exists an isotopy $\varphi : S^3 \times [0,1] \rightarrow S^3$ such that $\varphi_{0} = id_{S^3}$, $\varphi_1 (L_1) = L_2$ and $j_2 \circ \varphi_1 = \varphi_1 \circ j_1$. Here, $\varphi_{t} : S^3 \rightarrow S^3$ is defined as $\varphi_{t} (x) = \varphi (x,t)$.
Moreover, if $j_1 = j_2 = j$ and $j \circ \varphi_{t} = \varphi_{t} \circ j$ for every $t \in [0,1]$, then we call $(L_1, j)$ and $(L_2, j)$ are $j$-isotopic.
\end{definition}

\begin{remark}
The moves noted in Remark \ref{rem:moves} give $j$-isotopy for the associated links of divides with cusps.
\end{remark}

From now on, we fix the involution $j:S(D^2) \rightarrow S(D^2) ; j (x,v) = (x,-v)$ defined as before. 
The following proposition is clear by the definition of the associated link of a divide with cusps.

\begin{proposition}\label{prop:onlyif}
Let $P$ be a divide with cusps. The associated link $L(P)$ satisfies $j(L(P)) = L(P)$.
\end{proposition}

More precisely, we have the following proposition by decomposing a divide into components.
Let $P=P_1 \cup P_2 \cup \cdots \cup P_m$ be a component decomposition of a divide.

\begin{proposition}\label{prop:component}
Let $P_i$ be a component of a divide $P$. If $P_i$ is an interval component, then the link $L(P_i)$ is a strongly invertible knot. If $P_i$ is a circle component with odd cusps, then the link $L(P_i)$ is a $2$-periodic knot. 
If $P_i$ is a circle component with even cusps, then the link $L(P_i)$ is a $2$-component link and it is either a strongly invertible link or a $2$-periodic link.
\end{proposition}



\proof
We give proof only for the circle component case. 
The interval component case is proven similarly.
Let $P_i$ be a circle component with $\ell$ cusps. 
Let $S^1 =\{z \in \bC \mid |z| =1\}$ be a unit circle of a complex plane.
The component $P_i$ is described as the image of a continuous map $\overline{c}:S^1 \rightarrow P_i$ satisfying following conditions:
\begin{itemize}
\item The point $\overline{c}(1)$ is a smooth point.
\item Each cusp point $x_k$ is expressed as $x_k = \overline{c}(\zeta^k)$, where $\zeta = \exp (2\sqrt{-1} \pi/ (\ell + 1))$.
\item Define $c_{k} : A_{k} \rightarrow P_i$ as $c_{k} (z) = \overline{c} (z)$ for $1 \leq k \leq \ell+1$. Then, $c_{k}$ is an oriented smooth parametrization of a curve. Here $A_{k}$ is the arc $\{z \in S^1 \mid \frac{2(k-1)\pi}{(\ell + 1)} < \arg z < \frac{2k \pi }{ (\ell  +1) }\}$.
\end{itemize}
Let us define a partial parametrization of the link $\gamma^{\pm}_{k} : A_{k} \rightarrow L(P_i)$ for $1 \leq k \leq \ell +1$ as 
\[
\gamma^{\pm}_{k} (z) = (c_{k} (z) , \pm \frac{\sqrt{1-|c_{k}(z)|^2}}{|\dot c_{k}(z)|} \cdot \dot c_{k}(z)) .
\]
Remark that $j(\gamma^{\pm}_{k} (z)) = \gamma^{\mp}_{k}(z)$. At this time, we have the oriented smooth parametrization of the link $L(P_{i})$ as follows.
Note that we can extend the parametrization to cusps smoothly.

At first, suppose that the number of cusps $\ell$ is odd. Let us define $\Gamma : S^1 \rightarrow L(P_{i})$ as follows:
\begin{eqnarray*}
\Gamma (z) =  
\left \{
\begin{array}{ll}
\gamma^{+}_{k} (z^2) & (\mbox{$z \in A_{2k+1}$ and $0\leq \arg z \leq \pi$, $z \in A_{2k}$ and $\pi \leq \arg z \leq 2\pi$}), \\
\gamma^{-}_{k} (z^2) & (\mbox{$z \in A_{2k}$ and $0\leq \arg z \leq \pi$, $z \in A_{2k+1}$ and $\pi \leq \arg z \leq 2\pi$}). 
\end{array}
\right .
\end{eqnarray*}

This gives an oriented parametrization of $L(P_{i})$. 
By applying the involution, we have
\begin{eqnarray*}
j(\Gamma (z)) =  
\left \{
\begin{array}{ll}
\gamma^{-}_{k} (z^2) & (\mbox{$z \in A_{2k+1}$ and $0\leq \arg z \leq \pi$, $z \in A_{2k}$ and $\pi \leq \arg z \leq 2\pi$}), \\
\gamma^{+}_{k} (z^2) & (\mbox{$z \in A_{2k}$ and $0\leq \arg z \leq \pi$, $z \in A_{2k+1}$ and $\pi \leq \arg z \leq 2\pi$}). 
\end{array}
\right .
\end{eqnarray*}
thus $j(\Gamma(z)) = \Gamma (-z) $. Since the correspondence $z \mapsto -z$ preserves the orientation, the associated link $L(P_i)$ is $2$-periodic if the number of cusps is odd.

Next, suppose that the number of cusps $\ell$ is even. In this case, the link $L(P_i)$ consists of two components, and we denote by $L_1$ and $L_2$.
First, we give an oriented parametrization $\Gamma: S^1 \rightarrow L_1$ as follows.
\begin{eqnarray*}
\Gamma_1 (z) =  
\left \{
\begin{array}{ll}
\gamma^{+}_{k} (z) & (z \in A_{2k+1}) \\
\gamma^{-}_{k} (z) & (z \in A_{2k}) 
\end{array}
\right .
\end{eqnarray*}
Under fixing the orientation for $L_1$, we can choose two orientations for another component $L_{2}$.
First, define an oriented parametrization $\Gamma_2 : S^1 \rightarrow L_{2}$ as follows.
\begin{eqnarray*}
\Gamma_{2}(t) =
\left \{
\begin{array}{ll} 
\gamma^{-}_{k} (z) & (z \in A_{2k+1}) \\
\gamma^{+}_{k} (z) & (z \in A_{2k})
\end{array}
\right .
\end{eqnarray*}
Then, $j(\Gamma_{1} (z)) = \Gamma_{2} (z)$ and $j(\Gamma_{2} (z)) = \Gamma_{1} (z)$. Therefore, the orientation of the link $L(P_i)$ is preserved by the involution and thus it is a $2$-periodic link.

The other oriented parametrization $\Gamma_2 : S^1 \rightarrow L_{2}$ is defined as follows.
\begin{eqnarray*}
\Gamma_{2}(t) =
\left \{
\begin{array}{ll} 
\gamma^{-}_{k} (\overline{z}) & (z \in A_{2k+1}) \\
\gamma^{+}_{k} (\overline{z}) & (z \in A_{2k})
\end{array}
\right .
\end{eqnarray*}
Here, $\overline{z}$ is the complex conjugation of $z$. This time, $j(\Gamma_{1} (z)) = \Gamma_{2} (\overline{z})$ and $j(\Gamma_{2} (z)) = \Gamma_{1} (\overline{z})$. Since the correspondence $z \mapsto \overline{z}$ reverses the orientation, the link $L(P_{i})$ is a strongly invertible link.

\endproof

\begin{remark}\label{rem:compati}
As in the proof, from a circle component with even cusps, we can obtain either a $2$-periodic link or a strongly invertible link. This depends on just the choice of the orientation of the associated link.
\end{remark}


\section{Main result and examples}
\label{sec:main}

In this section, we give the proof of the main theorem. The main theorem is the following.

\begin{theorem}\label{thm:main}
A link $L$ in $S^3$ is described as the associated link of a divide with cusps if and only if there exists an involution $j$ on $S^3$ which has a non-empty fixed point set such that $j(L) = L$.
\end{theorem}

The "only if" part is nothing but Proposition \ref{prop:onlyif}. Thus, it suffices to prove the following theorem.  

\begin{theorem}
Suppose that link $L$ satisfies $j(L) = L$ for an involution $j$ on $S^3$ with a non-empty fixed point set. Then, there exists a divide with cusps $P$ such that $L(P)$ is $j$-isotopic to $L$.
\end{theorem}

\proof
Let $L$ be a link in $S^3$ identified with the sphere $S(D^2)$. 
Suppose that $L$ satisfies $j(L) = L$ for the involution $j(x,v) = (x,-v)$.
Let $P = \pi(L)$ be the image by the projection.
We can generically suppose that $P$ satisfies the following conditions.
\begin{itemize}
\item [(i)] The image $P$ is a divide with no cusps, i.e., the union of the proper generic immersed curves.
\item [(ii)] For each $x \in P \setminus \partial P$, $L$ intersects to $\pi^{-1} (x)$ transversally.
\item[(iii)] The restriction $\pi |_{L\cap \pi^{-1} (P_{0})} : L \cap \pi^{-1}(P_0) \rightarrow P_{0} $ is a branched $2$-fold covering ramified at $\partial P_0$, where $P_0$ is obtained from $P$ by removing its cross points.
\end{itemize}
Since $\Fix (j) = \partial D^2 \times \{0\}$, a fixed point in the link $L$ by the involution can be written as $(e^{\sqrt{-1} \theta},0) \in \Fix (j) \cap L$.
The unit tangent vector of $L$ at this point can be expressed as $(0,e^{\sqrt{-1} \psi})$. By rotating around $\Fix (j)$, we can assume that $\psi = \theta + \pi$.
Therefore, by performing $j$-isotopy, we have that the link $L$ coincides to $L(P)$ at the $\delta$ neighborhood $\{p \in D^2 \mid |p| \geq 1 -\delta \}$ of $\partial D^2$ for a sufficiently small value $\delta > 0$.

\begin{figure}[htbp]
\centering 
\begin{tikzpicture}
\coordinate (P) at (0,0);
\coordinate (Q) at (-5,0);
\coordinate (R) at (7,0);


\draw (Q) circle (1);

\draw[->,>=stealth, thick, red] (Q) --++(0.945,0.315);
\draw[->,>=stealth, thick, red, densely dotted] (Q) --++(-0.945,-0.315);
\draw[->,>=stealth, thick, blue] (Q) --++(-0.315,0.945);
\draw[->,>=stealth, thick, blue, densely dotted] (Q) --++(0.315,-0.945);
\fill (Q) circle (0.06);

\draw (Q) ++ (0.945,0.315) node[right]{$\theta_{1}$};
\draw (Q) ++ (-0.315,0.945) node[above]{$\theta_{2}$};

\draw[->,>=stealth] (Q) --++ (1,0) ;
\draw (Q)++(1,-0.1)node[right]{$0$};

\draw[thick] (P) circle (2);

\draw (P) ++(1.41,1.41) node[above]{$a_{1}$};
\draw (P) ++(-1.41,1.41) node[left]{$a_{2}$};
\draw (P) ++(-1.41,-1.41) node[below]{$a_{3}$};
\draw (P) ++(1.41,-1.41) node[right]{$a_{4}$};

\draw (P) ++(1.41,1.41) --++(-2.82,-2.82);
\draw (P) ++(-1.41,1.41) --++(2.82,-2.82);

\draw[->,>=stealth, thick, red] (P) --++(0.3,0.1);
\draw[->,>=stealth, thick, red, densely dotted] (P) --++(-0.3,-0.1);
\draw[->,>=stealth, thick, red] (P)++(1.41,1.41) --++(0.3,0.1);
\draw[->,>=stealth, thick, red, densely dotted] (P)++(1.41,1.41) --++(-0.3,-0.1);
\draw[->,>=stealth, thick, red] (P)++(-1.41,-1.41) --++(0.3,0.1);
\draw[->,>=stealth, thick, red, densely dotted] (P)++(-1.41,-1.41) --++(-0.3,-0.1);
\draw[->,>=stealth, thick, red] (P)++(0.7,0.7) --++(0.3,0.1);
\draw[->,>=stealth, thick, red, densely dotted] (P)++(0.7,0.7) --++(-0.3,-0.1);
\draw[->,>=stealth, thick, red] (P)++(-0.7,-0.7) --++(0.3,0.1);
\draw[->,>=stealth, thick, red, densely dotted] (P)++(-0.7,-0.7) --++(-0.3,-0.1);

\draw[->,>=stealth, thick, blue] (P) --++(-0.1,0.3);
\draw[->,>=stealth, thick, blue, densely dotted] (P) --++(0.1,-0.3);
\draw[->,>=stealth, thick, blue] (P) ++(-1.41,1.41) --++(-0.1,0.3);
\draw[->,>=stealth, thick, blue, densely dotted] (P)++(-1.41,1.41) --++(0.1,-0.3);
\draw[->,>=stealth, thick, blue] (P)++(1.41,-1.41) --++(-0.1,0.3);
\draw[->,>=stealth, thick, blue, densely dotted] (P)++(1.41,-1.41) --++(0.1,-0.3);
\draw[->,>=stealth, thick, blue] (P)++(-0.7,0.7) --++(-0.1,0.3);
\draw[->,>=stealth, thick, blue, densely dotted] (P)++(-0.7,0.7)  --++(0.1,-0.3);
\draw[->,>=stealth, thick, blue] (P)++(0.7,-0.7)  --++(-0.1,0.3);
\draw[->,>=stealth, thick, blue, densely dotted] (P)++(0.7,-0.7)  --++(0.1,-0.3);

\draw (P) ++ (0,-0.2) node[below] {$x$};

\draw [thick](R) circle (2);

\draw (R) ++(1.41,1.41) node[above]{$a_{1}$};
\draw (R) ++(-1.41,1.41) node[left]{$a_{2}$};
\draw (R) ++(-1.41,-1.41) node[below]{$a_{3}$};
\draw (R) ++(1.41,-1.41) node[right]{$a_{4}$};

\draw (R) ++(1.41,1.41) --++(-2.82,-2.82);
\draw (R) ++(-1.41,1.41) --++(2.82,-2.82);

\draw[->,>=stealth, thick, red] (R)++(0.02,-0.02) --++(0.223,0.223);
\draw[->,>=stealth, thick, red, densely dotted] (R)++(-0.02,0.02) --++(-0.223,-0.223);
\draw[->,>=stealth, thick, red] (R)++(0.02,-0.02) ++(1.41,1.41) --++(0.223,0.223);
\draw[->,>=stealth, thick, red, densely dotted] (R)++(-0.02,0.02) ++(1.41,1.41) --++(-0.223,-0.223);
\draw[->,>=stealth, thick, red] (R)++(0.02,-0.02) ++(-1.41,-1.41) --++(0.223,0.223);
\draw[->,>=stealth, thick, red, densely dotted] (R)++(-0.02,0.02) ++(-1.41,-1.41) --++(-0.223,-0.223);
\draw[->,>=stealth, thick, red] (R)++(0.02,-0.02) ++(0.7,0.7) --++(0.223,0.223);
\draw[->,>=stealth, thick, red, densely dotted] (R)++(-0.02,0.02) ++(0.7,0.7) --++(-0.223,-0.223);
\draw[->,>=stealth, thick, red] (R)++(0.02,-0.02) ++(-0.7,-0.7) --++(0.223,0.223);
\draw[->,>=stealth, thick, red, densely dotted] (R)++(-0.02,0.02) ++(-0.7,-0.7) --++(-0.223,-0.223);

\draw[->,>=stealth, thick, blue] (R)++(0.02,0.02) --++(-0.223,0.223);
\draw[->,>=stealth, thick, blue, densely dotted] (R)++(-0.02,-0.02) --++(0.223,-0.223);

\draw[->,>=stealth, thick, blue] (R)++(0.02,0.02)++(-1.41,1.41) --++(-0.223,0.223);
\draw[->,>=stealth, thick, blue, densely dotted] (R)++(-0.02,-0.02)++(-1.41,1.41) --++(0.223,-0.223);
\draw[->,>=stealth, thick, blue] (R)++(0.02,0.02)++(1.41,-1.41) --++(-0.223,0.223);
\draw[->,>=stealth, thick, blue, densely dotted] (R)++(-0.02,-0.02)++(1.41,-1.41) --++(0.223,-0.223);
\draw[->,>=stealth, thick, blue] (R)++(0.02,0.02)++(-0.7,0.7) --++(-0.223,0.223);
\draw[->,>=stealth, thick, blue, densely dotted] (R)++(-0.02,-0.02)++(-0.7,0.7) --++(0.223,-0.223);
\draw[->,>=stealth, thick, blue] (R)++(0.02,0.02)++(0.7,-0.7) --++(-0.223,0.223);
\draw[->,>=stealth, thick, blue, densely dotted] (R)++(-0.02,-0.02)++(0.7,-0.7) --++(0.223,-0.223);

\draw[->,>=stealth, thick, blue] (P) ++(-1.41,1.41) --++(-0.1,0.3);
\draw[->,>=stealth, thick, blue, densely dotted] (P)++(-1.41,1.41) --++(0.1,-0.3);
\draw[->,>=stealth, thick, blue] (P)++(1.41,-1.41) --++(-0.1,0.3);
\draw[->,>=stealth, thick, blue, densely dotted] (P)++(1.41,-1.41) --++(0.1,-0.3);
\draw[->,>=stealth, thick, blue] (P)++(-0.7,0.7) --++(-0.1,0.3);
\draw[->,>=stealth, thick, blue, densely dotted] (P)++(-0.7,0.7)  --++(0.1,-0.3);
\draw[->,>=stealth, thick, blue] (P)++(0.7,-0.7)  --++(-0.1,0.3);
\draw[->,>=stealth, thick, blue, densely dotted] (P)++(0.7,-0.7)  --++(0.1,-0.3);

\draw (R) ++ (0,-0.2) node[below] {$x$};

\draw[->] (P) ++ (2.5,0) --++(2,0);
\draw (P) ++ (3.5,0) node[below] {move by};
\draw (P) ++ (3.5,-0.4) node[below] {$j$-isotopy $\varphi$};

\end{tikzpicture}

\label{fig:crossing}
\caption{Around the double point}
\end{figure}
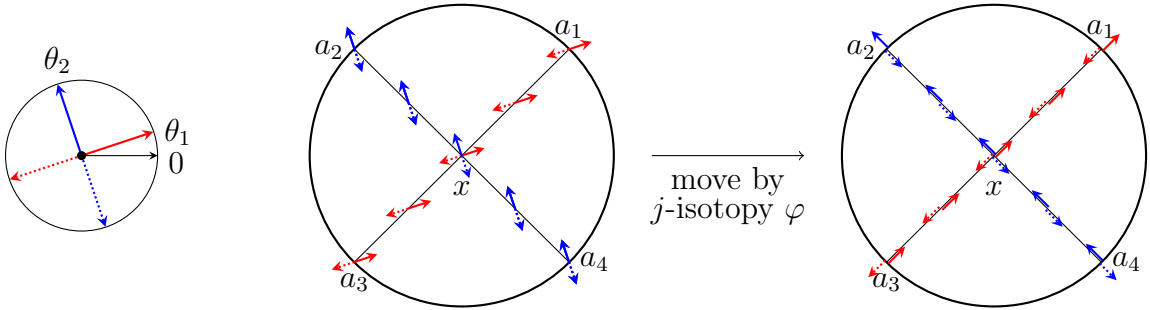

Next, we consider the neighborhood of the double point of $P$. Let $x \in P$ be a double point and $D_{\delta}$ be a disk centered at $x$ with the radius $\delta$.
Let $a_1, a_2, a_3$ and $a_4$ be the intersection of $P$ and $\partial D_{\delta}$. 
We suppose that $a_{i}$'s are ordered counterclockwisely and $0 \leq \alpha_1 < \alpha_2 < \alpha_3 < \alpha_4 < 2 \pi$, where $\alpha_{i} = \arg (a_{i}-x)$.
We can assume that the intersection $P \cap D_{\delta}$ is the union of intervals $[a_1,a_3]$ and $[a_2,a_4]$.
The link $L \cap \pi^{-1} (D_{\delta})$ is expressed as $L \cap \pi^{-1} ([a_1,a_3] \cup [a_2, a_4])$, and by moving $j$-isotopy, each component is described as follows ($i=1,2$).
\[
L \cap \pi^{-1} ([a_i, a_{i+2}]) = \{(x,\pm \theta_{i}) \mid x \in [a_{i},a_{i+2}]\},
\]
where $0 \leq \theta_{1} < \theta_{2} < \pi$.
The link $L \cap \pi^{-1} (D_{\delta}) $ consists of four segments and each component can be parametrized as follows:
\[
\gamma_{i} (s) = (s a_{i+2} + (1-s) a_{i} , \varepsilon \theta_{i})
\]
where $0 \leq s \leq 1$ , $a_5 = a_1$ , $a_6 = a_2$. 
The signature $\varepsilon = 1$ if $i=1,2$ and $\varepsilon = -1$ if $i=3,4$. Note that this parametrization may not agree with the original orientation of the link.
Then, we can define a $j$-isotopy $\{ \varphi_{t} \}_{t \in [0,1]}$ as follows.
\[
\varphi_{t} (\gamma_i (s)) = \left ( (s a_{i+2} + (1-s) a_{i}) , (1-t) \varepsilon a_{i} + t \alpha_{i}  \right )
\]
It folows that $\varphi_{0} (\gamma_{i}(s)) = \gamma_{i} (s)$ and $\varphi_{1} (\gamma_{i} (s)) = (s a_{i+2} + (1-s) a_{i}, \alpha_{i})$. Therefore, the link $L \cap \pi^{-1} (D_{\delta})$ coincides with the divide link $L(P)$ up to $j$-isotopy (see Figure \ref{fig:crossing}). 


Let $P'$ be the set obtained from $P$ by removing all the $\delta$-neighborhoods of crossing points and $\delta$-neighborhood of $\partial D^2$. 
The set $P'$ is the union of finite embedded intervals and embedded circles.
Let $E$ be a connected component of $P'$ and suppose that $E$ is an interval.
The link $L \cap \pi^{-1} (E)$ lives in the annulus $\pi^{-1} (E) \cong E \times S^1$ and $L \cap \pi^{-1} (E)$ coincides with $L(P')$ at the boundary $L \cap \pi^{-1} (\partial E)$.

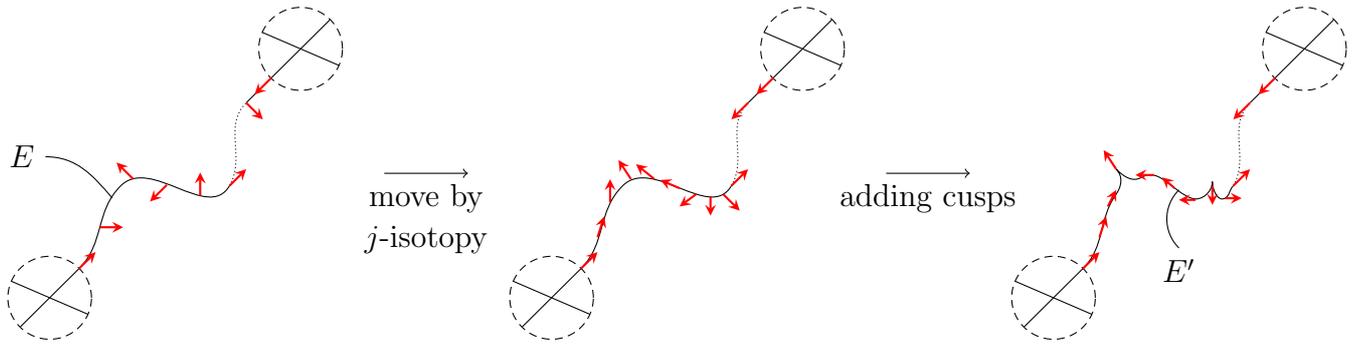
\begin{figure}
\centering
\begin{tikzpicture}[scale =1.1]
\coordinate (P1) at (0,0);
\coordinate (P2) at (3,3);
\coordinate (Q1) at (6,0);
\coordinate (Q2) at (9,3);
\coordinate (R1) at (12,0);
\coordinate (R2) at (15,3);

\draw[densely dashed] (P1) circle (0.5);
\draw[densely dashed] (P2) circle (0.5);

\draw (P1) --++(0.46,-0.2);
\draw (P1) --++(-0.46,0.2);
\draw (P1) --++(0.35,0.35);
\draw (P1) --++(-0.35,-0.35);

\draw (P2) --++(0.46,-0.2);
\draw (P2) --++(-0.46,0.2);
\draw (P2) --++(0.35,0.35);
\draw (P2) --++(-0.35,-0.35);

\draw (P1) ++ (0.35,0.35) to[out=45,in=-135] ++(0.5,1) to[out=45,in=-120] ++(1.3,0)to[out=60,in=-135] ++(0.2,1)--++(0.3,0.3);
\draw[preaction={draw = white, line width = 1 pt},densely dotted]  (P1) ++ (0.35,0.35)++(0.5,1) ++(1.3,0) to[out=60,in=-135] ++(0.2,1);

\draw (P1) ++ (0.75,1.2) to[out=130,in=0] ++(-0.8,0.5) node[left]{$E$};

\draw[->,>=stealth, red, thick] (P1) ++ (0.35,0.35) --++(0.2,0.2);
\draw[->,>=stealth, red, thick] (P1) ++ (0.6,0.85) --++(0.28,0);
\draw[->,>=stealth, red, thick] (P1) ++ (1,1.43) --++(-0.2,0.2);
\draw[->,>=stealth, red, thick] (P1) ++ (1.4,1.36) --++(-0.2,-0.2);
\draw[->,>=stealth, red, thick] (P1) ++ (1.8,1.23) --++(0,0.28);
\draw[->,>=stealth, red, thick] (P1) ++ (2.15,1.35) --++(0.2,0.2);

\draw[->,>=stealth, red, thick] (P1) ++ (2.35,2.35) --++(0.2,-0.2);
\draw[->,>=stealth, red, thick] (P2) ++ (-0.35,-0.35) --++(-0.2,-0.2);


\draw[densely dashed] (Q1) circle (0.5);
\draw[densely dashed] (Q2) circle (0.5);

\draw (Q1) --++(0.46,-0.2);
\draw (Q1) --++(-0.46,0.2);
\draw (Q1) --++(0.35,0.35);
\draw (Q1) --++(-0.35,-0.35);

\draw (Q2) --++(0.46,-0.2);
\draw (Q2) --++(-0.46,0.2);
\draw (Q2) --++(0.35,0.35);
\draw (Q2) --++(-0.35,-0.35);

\draw (Q1) ++ (0.35,0.35) to[out=45,in=-135] ++(0.5,1) to[out=45,in=-120] ++(1.3,0)to[out=60,in=-135] ++(0.2,1)--++(0.3,0.3);
\draw[preaction={draw = white, line width = 1 pt},densely dotted]  (Q1) ++ (0.35,0.35)++(0.5,1) ++(1.3,0) to[out=60,in=-135] ++(0.2,1);

\draw[->,>=stealth, red, thick] (Q1) ++ (0.35,0.35) --++(0.2,0.2);
\draw[->,>=stealth, red, thick] (Q1) ++ (0.55,0.73) --++(0.08,0.25);
\draw[->,>=stealth, red, thick] (Q1) ++ (0.7,1.15) --++(0,0.28);
\draw[->,>=stealth, red, thick] (Q1) ++ (0.95,1.43) --++(-0.15,0.23);
\draw[->,>=stealth, red, thick] (Q1) ++ (1.23,1.43) --++(-0.22,0.18);
\draw[->,>=stealth, red, thick] (Q1) ++ (1.53,1.32) --++(-0.23,0.1);

\draw[->,>=stealth, red, thick] (Q1) ++ (1.73,1.25) --++(-0.2,-0.15);
\draw[->,>=stealth, red, thick] (Q1) ++ (1.9,1.22) --++(0,-0.23);
\draw[->,>=stealth, red, thick] (Q1) ++ (2.05,1.25) --++(0.2,-0.2);

\draw[->,>=stealth, red, thick] (Q1) ++ (2.15,1.35) --++(0.2,0.2);

\draw[->,>=stealth, red, thick] (Q1) ++ (2.35,2.35) --++(-0.2,-0.2);
\draw[->,>=stealth, red, thick] (Q2) ++ (-0.35,-0.35) --++(-0.2,-0.2);


\draw[densely dashed] (R1) circle (0.5);
\draw[densely dashed] (R2) circle (0.5);

\draw (R1) --++(0.46,-0.2);
\draw (R1) --++(-0.46,0.2);
\draw (R1) --++(0.35,0.35);
\draw (R1) --++(-0.35,-0.35);

\draw (R2) --++(0.46,-0.2);
\draw (R2) --++(-0.46,0.2);
\draw (R2) --++(0.35,0.35);
\draw (R2) --++(-0.35,-0.35);


\draw (R2) ++ (-0.35,-0.35) --++(-0.3,-0.3);

\draw (R1) ++ (0.35,0.35) to[out=45,in=-125] ++(0.4,0.9)to[out=55,in=-45] ++(0,0.3)to[out=-45,in=-150]++(0.3,-0.1)to[out=30,in=160]  ++(0.6,-0.25)to[out=-20,in=-90]++(0.25,0.2)to[out=-90,in=-150] ++(0.15,-0.2)to[out=30,in=-135]++(0.1,0.15);

\draw[preaction={draw = white, line width = 1 pt},densely dotted]  (R1) ++ (0.35,0.35)++(0.5,1) ++(1.3,0) to[out=60,in=-135] ++(0.2,1);

\draw[->,>=stealth, red, thick] (R1) ++ (0.35,0.35) --++(0.2,0.2);
\draw[->,>=stealth, red, thick] (R1) ++ (0.55,0.73) --++(0.08,0.25);
\draw[->,>=stealth, red, thick] (R1) ++ (0.65,1.1) --++(0.1,0.2);
\draw[->,>=stealth, red, thick] (R1) ++ (0.75,1.55) --++(-0.15,0.23);
\draw[->,>=stealth, red, thick] (R1) ++ (1.2,1.48) --++(-0.22,-0);
\draw[->,>=stealth, red, thick] (R1) ++ (1.45,1.32) --++(-0.15,0.13);

\draw[->,>=stealth, red, thick] (R1) ++ (1.7,1.18) --++(-0.2,-0);
\draw[->,>=stealth, red, thick] (R1) ++ (1.9,1.35) --++(0,-0.23);
\draw[->,>=stealth, red, thick] (R1) ++ (2.05,1.2) --++(0.2,-0);

\draw[->,>=stealth, red, thick] (R1) ++ (2.15,1.35) --++(0.2,0.2);

\draw[->,>=stealth, red, thick] (R1) ++ (2.35,2.35) --++(-0.2,-0.2);
\draw[->,>=stealth, red, thick] (R2) ++ (-0.35,-0.35) --++(-0.2,-0.2);

\draw (R1) ++(1.5,1.3) to[out=-135,in=135]++(0,-0.7)node[below]{$E'$};

\draw[->] (P2) ++(1,-1.5) --++(1,0); 
\draw (P2)++(1.5,-1.5) node[below]{move by} ;
\draw (P2)++(1.5,-2) node[below]{$j$-isotopy} ;

\draw[->] (Q2) ++(1,-1.5) --++(1,0); 
\draw (Q2)++(1.5,-1.5) node[below]{adding cusps} ;

\end{tikzpicture}

\caption{The interval component}
\label{fig:interval}
\end{figure}

Let $\Phi$ be a self homeomorphism of $\pi^{-1} (E)$ defined by a half twist of $\pi^{-1} (E)$. 
Then, there exists an integer $n$ such that $L \cap \pi^{-1} (E)$ is equivalent to $\Phi^{n} (L(P) \cap \pi^{-1} (E))$. Bringing half twists to the link of a divide corresponds to adding cusps to the divide as noted in Remark \ref{rem:halftwists}. 
Let $E'$ be the interval with the corresponding cusp added.
Then, the link $L \cap \pi^{-1} (E)$ coincides with the associated link of $L(E')$. (see Figure \ref{fig:interval})

When $E$ is an embedded circle component, by a similar argument, it follows that the link $L \cap \pi^{-1} (E)$ coincides with the associated link of a divide with finite cusps added.
Let $\overline{P}$ be the resulting divide with cusps, obtained from $P$ by adding finite cusps according to all the previous procedures. 
Then, the link $L$ coincides with the associated link $L(\overline{P})$.
\endproof

As the special cases of the result, we have the following corollaries by combining Proposition \ref{prop:component} and Theorem \ref{thm:main}. 

\begin{corollary}
Every strongly invertible link is described as the associated link of a divide with cusps. Particularly, it is the associated link of a divide with cusps consisting of interval components and circle components with even cusps. 
\end{corollary}

\begin{corollary}
Every $2$-periodic link is described as the associated link of a divide with cusps. Particularly, it is the associated link of a divide with cusps consisting of circle components. 
\end{corollary}


Let us consider the case where $L=K$ is a knot in $S^3$. If a knot $K$ satisfies $j(K) = K$, $K$ is either a strongly invertible knot or a $2$-periodic knot. 
Thus, we have the following theorem as the specialization of the main result.

\begin{theorem}
A knot $K$ is described as the associated link of a divide with cusps if and only if $K$ is either a strongly invertible knot or $2$-periodic knot.
\end{theorem}

\begin{figure}[htbp]
\centering

\begin{tikzpicture}
\coordinate (P) at (0,0);
\coordinate (Q) at (-1,-7);
\coordinate (R) at (3.5,-7);
\coordinate (S) at (7.5,-7);

\coordinate (P1) at (6,-3);
\coordinate (P2) at (6,-1);
\coordinate (P3) at (6,1);
\coordinate (P4) at (6,3);

\draw[thick] (P) circle (1.5);
\draw (P) ++(-1.05,-1.05) --++ (2.1,2.1);
\draw (P) ++(-1.05,1.05) --++ (2.1,-2.1);

\draw[->,>=stealth] (P) --++(0.1,0.4)node[above]{$d$};
\draw[->,>=stealth] (P) --++(-0.1,-0.4);

\draw (P) node[right]{$-$};

\draw (P) ++ (0,-2) node{(1)};


\draw[<-,>=stealth] (Q) ++ (-1.2,-0.8) --++(2,2);
\draw[draw=white, line width = 5 pt] (Q) ++ (-0.8,1.2) --++(2,-2);
\draw[<-,>=stealth] (Q) ++ (-0.8,1.2) --++(2,-2);
\draw[draw=white, line width = 5 pt] (Q) ++ (-0.8,-1.2) --++(2,2);
\draw[->,>=stealth] (Q) ++ (-0.8,-1.2) --++(2,2);
\draw[draw=white, line width = 5 pt] (Q) ++ (-1.2,0.8) --++(2,-2);
\draw[->,>=stealth] (Q) ++ (-1.2,0.8) --++(2,-2);

\draw (Q) ++ (2,-2) node{(3)};

\draw[->] (Q) ++ (1.5,0) --++(1.5,0) ;
\draw (Q) ++(2.25,0) node[below]{move by};
\draw (Q) ++(2.25,-0.4) node[below]{isotopy};


\draw [<-,>=stealth] (R) ++ (-1.2,-0.8) --++(2,2); 

\draw[draw=white, line width = 5 pt] (R) ++ (-0.8,1.2) --++(0.1,-0.1) to[out=-45,in=135] ++ (0,-0.8) --++(1,-1) to[out=-45,in=135] ++ (0.8,0)--++(0.1,-0.1);
\draw[<-,>=stealth] (R) ++ (-0.8,1.2) --++(0.1,-0.1) to[out=-45,in=135] ++ (0,-0.8) --++(1,-1) to[out=-45,in=135] ++ (0.8,0)--++(0.1,-0.1);

\draw[draw=white, line width = 5 pt] (R) ++ (-0.8,-1.2) --++(2,2);
\draw[->,>=stealth] (R) ++ (-0.8,-1.2) --++(2,2);

\draw[draw=white, line width = 5 pt] (R) ++ (-1.2,0.8) --++(0.1,-0.1)to[out=-45,in=135] ++ (0.8,-0) --++(1,-1) to[out=-45,in=135] ++ (0,-0.8)--++(0.1,-0.1) ;
\draw[->,>=stealth] (R) ++ (-1.2,0.8) --++(0.1,-0.1)to[out=-45,in=135] ++ (0.8,-0) --++(1,-1) to[out=-45,in=135] ++ (0,-0.8)--++(0.1,-0.1) ;

\draw [->,>=stealth] (R) ++ (-1.2,0.8) --++(0.1,-0.1) ;
\draw [->,>=stealth] (R) ++ (1.2,-0.8) --++(-0.1,0.1) ;
\draw [->,>=stealth] (R) ++ (1.2,-0.8) ++(-1.4,1.4) --++(-0.1,0.1) ;
\draw [->,>=stealth] (R) ++ (-1.2,0.8) ++(1.4,-1.4) --++(0.1,-0.1) ;

\draw[thick] (S) circle (1.5);
\draw (S) ++(-1.05,-1.05) --++ (2.1,2.1);
\draw (S) ++(-1.05,1.05) --++ (0.3,-0.3)to[out=-45,in=-135] ++(0.5,0)to[out=-135,in=135]++(0,-0.5) --++(0.5,-0.5)to[out=-45,in=45] ++(0,-0.5)to[out=45,in=135]++(0.5,-0)--++(0.3,-0.3);

\draw (S) ++ (0,-2) node{(4)};




\draw (P1)circle [x radius=2,y radius=0.6];
\draw (P1) ++(-1.15,-0.5) --++(2.2,1);
\draw (P1) ++(-1.15,0.5) --++(2.2,-1);
\draw (P1) ++ (-2,0) node[left]{$-\frac{3}{4} \pi$};

\draw (P2)circle [x radius=2,y radius=0.6];
\draw (P2) ++(-1.15,-0.5) --++(2.2,1);
\draw (P2) ++(-1.15,0.5) --++(2.2,-1);
\draw (P2) ++ (-2,0) node[left]{$-\frac{1}{4} \pi$};

\draw (P3)circle [x radius=2,y radius=0.6];
\draw (P3) ++(-1.15,-0.5) --++(2.2,1);
\draw (P3) ++(-1.15,0.5) --++(2.2,-1);
\draw (P3) ++ (-2,0) node[left]{$\frac{1}{4} \pi$};

\draw (P4)circle [x radius=2,y radius=0.6];
\draw (P4) ++(-1.15,-0.5) --++(2.2,1);
\draw (P4) ++(-1.15,0.5) --++(2.2,-1);
\draw (P4) ++ (-2,0) node[left]{$\frac{3}{4} \pi$};

\draw (P4) ++ (-2,0.8)node[left]{$\pi$};
\draw (P1) ++ (-2,-0.8)node[left]{$-\pi$};

\draw [very thick] (P4) ++ (-0.4,0.7) to[out=-60,in=-150] ++ (0.3,-1) to[out=30, in=-120] ++ (0.4,1.2);
\draw [very thick] (P4) ++(-1.15,0.5) --++(0.77,-0.35);
\draw [preaction={draw=white, line width = 3 pt}, very thick] (P4) ++(-1.15,0.5) ++(0.77,-0.35) ++(0.22,-0.1) --++(1.21,-0.55);

\draw[very thick] (P3) ++(-1.15,-0.5) --++(0.55,0.25) to[out=20,in=100] ++ (0.4,-0.2) to[out=-80,in=90] ++ (-0.1,-1.7) to[out=-90,in=180] ++(0.1,-0.1) to[out=0,in=-120] ++(0.2,0.1);
\draw[very thick] (P3) ++(-1.15,-0.5) ++(2.2,1) --++(-0.55,-0.25) to[out=-160,in=80] ++(-0.2,-0.5) to[out=-100,in=80] ++(-0.23,-1.7);

\draw [very thick] (P2) ++(-1.15,0.5) --++(0.77,-0.35);
\draw [preaction={draw=white, line width = 3 pt}, very thick] (P2) ++(-1.15,0.5) ++(0.77,-0.35) ++(0.22,-0.1) --++(1.21,-0.55);

\draw[very thick] (P1) ++(-1.15,-0.5) --++(0.55,0.25) to[out=20,in=100] ++ (0.4,-0.2);
\draw[very thick] (P1) ++(-1.15,-0.5) ++(2.2,1) --++(-0.55,-0.25) to[out=-170,in=80] ++(-0.2,-0.5);

\draw (P1) ++(0,-1) node{$(2)$};

\draw[->, very thick] (P) ++(2,0) --++(1,0);
\draw[->, very thick] (P) ++(3,-4.5) --++(-1,-0.5);

\end{tikzpicture}
\caption{A negative sign of a signed divide and cusps in a divide with cusps}
\label{fig:negative}
\end{figure}
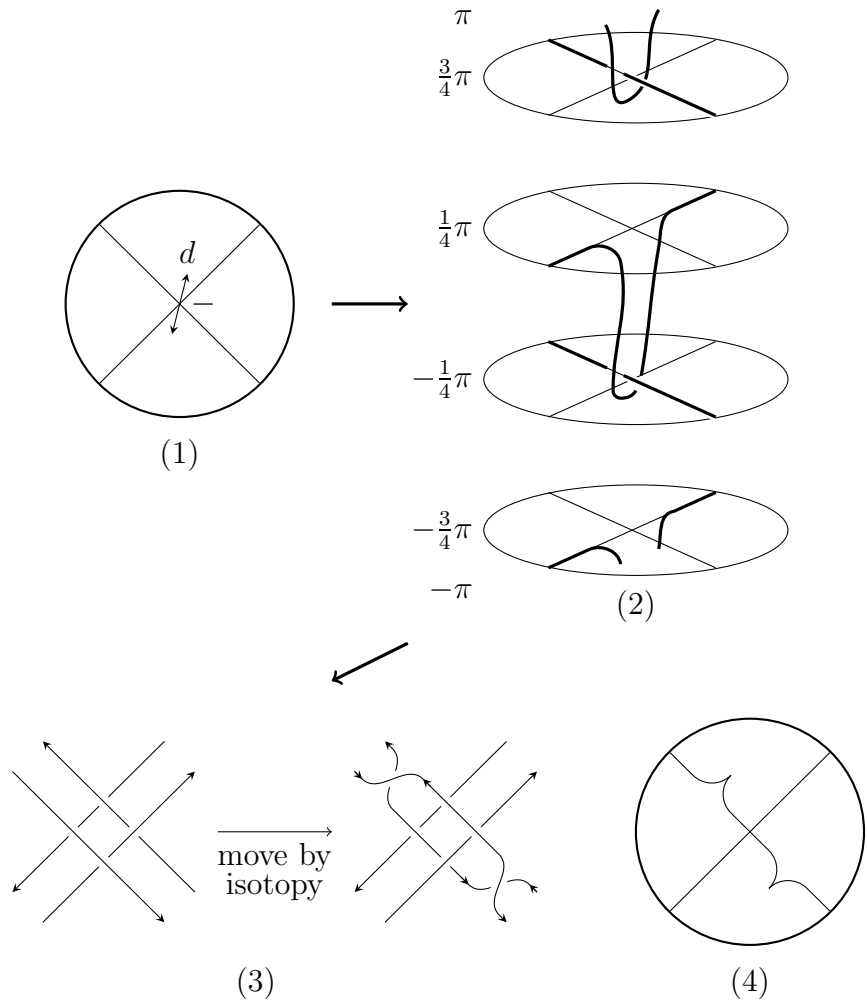

\begin{remark}\label{rem:negative}
A signed divide, which is introduced by Couture \cite{cou}, can be considered as the special case of the divide with cusps.
In fact, equipping a negative sign to a crossing point corresponds to adding some cusps. 
The rule of adding cusps is as in Figure \ref{fig:negative}.
If we give a negative sign to a double point as in $(1)$ in Figure \ref{fig:negative}, then the diagram of the associated link is locally drawn as in $(2)$, applying Hirasawa's algorithm. It can be projected to the diagram in the plane as in $(3)$ by a small perturbation. Here, we rotate the angle of the cylinder so that the branch of $- \pi /4 $ is at the top.
Moving by isotopy, it can be locally described as the associated link of a divide with cusps, by adding cusps appropriately around the double point as $(4)$.

Couture proved that every strongly invertible link is described as the associated link of a signed divide (\cite{cou}, Proposition 1.10).
Therefore, our main result is the generalization of this proposition.
\end{remark}

\begin{example}\label{exam:torus}
For $n \geq 1$, let $T(2,2n+1)$ be a $(2,2n+1)$-torus knot. There are two involutions $j_1, j_2$ of $S^3$ such that $(T(2,2n+1), j_1)$ is strongly invertible and $(T(2,2n+1), j_2)$ is $2$-periodic. 
Both links are described as the associated link of divides with cusps as in Figure \ref{fig:torusknot}.

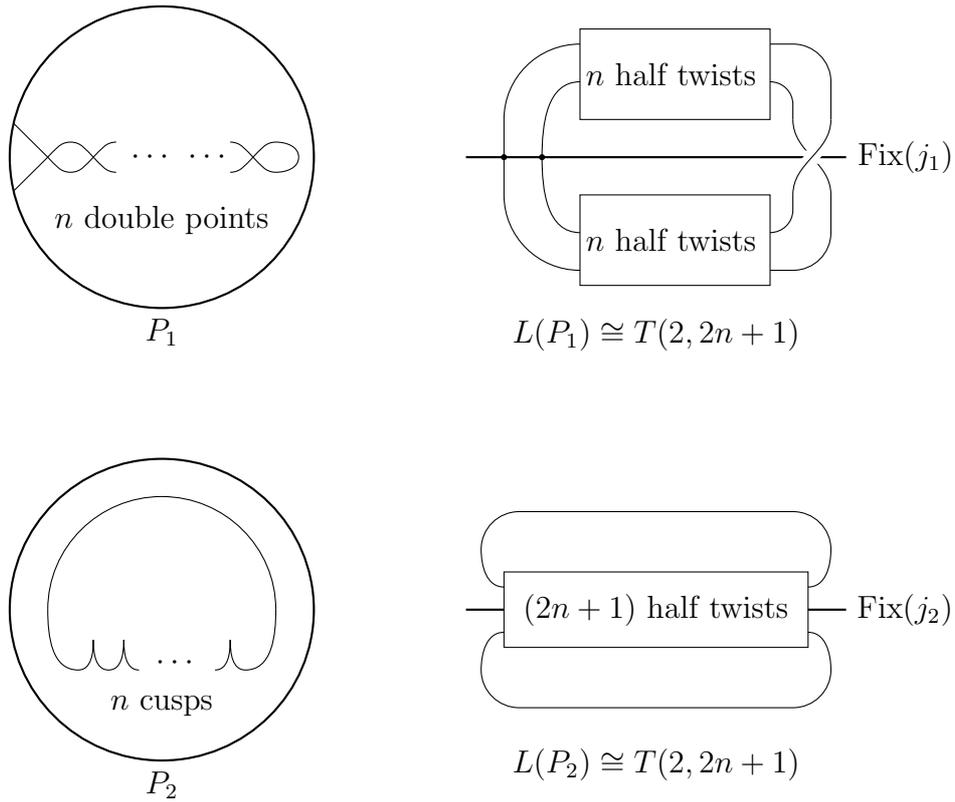
\begin{figure}[htbp]
\centering

\begin{tikzpicture}
\coordinate (P) at (0,0);
\coordinate (P1) at (4,0);
\coordinate (Q) at (0,-6);
\coordinate (Q1) at (4,-6);


\draw[thick] (P) circle (2);

\draw (P) ++(-1.5,0) --++(-0.45,0.45);
\draw (P) ++(-1.5,0) --++(-0.45,-0.45);

\draw (P) ++(-1.5,0) to[out=-45,in=180] ++(0.3,-0.2)to[out=0,in=-135]++(0.3,0.2)to[out=45,in=-180] ++(0.3,0.2);
\draw (P) ++(-1.5,0) to[out=45,in=180] ++(0.3,0.2)to[out=0,in=135]++(0.3,-0.2) to[out=-45,in=180] ++(0.3,-0.2);

\draw (P)++(0.3,0) node[left]{$\cdots$};
\draw (P)++(0.2,0) node[right]{$\cdots$};

\draw (P)++(1.2,0) to[out=45,in=180] ++ (0.3,0.2)to[out=0,in=90] ++(0.3,-0.2);
\draw (P)++(1.2,0) to[out=-45,in=180] ++ (0.3,-0.2)to[out=0,in=-90] ++(0.3,0.2);

\draw (P) ++ (1.2,0) to[out=135,in=0] ++(-0.3,0.2);
\draw (P) ++ (1.2,0) to[out=-135,in=0] ++(-0.3,-0.2);

\draw (P) ++ (0,-0.5) node[below]{$n$ double points};

\draw (P) ++(0,-2) node[below]{$P_1$};


\draw (P1) ++ (4.8,-0.5) to[out=90,in=-90] ++ (-0.5,1);

\draw[preaction={draw = white, line width =5pt}, thick] (P1) --++(5,0)node[right]{$\Fix(j_{1})$};

\fill(P1) ++(0.5,0) circle (0.04);
\fill(P1) ++(1,0) circle (0.04);

\draw (P1) ++ (0.5,0) --++(0,0.5)to[out=90,in=180] ++(1,1);
\draw (P1) ++ (1,0) to[out=90,in=180] ++(0.5,1);

\draw (P1) ++ (1.5,1.7) --++(2.5,0) --++(0,-1.2) --++(-2.5,0) --++(0,1.2)--cycle;

\draw (P1) ++ (4,1.5)--++(0.3,0)  to[out=0,in=90] ++ (0.5,-0.5)--++(0,-0.5);
\draw (P1) ++(4,1) to[out=0,in=90] ++(0.3,-0.3)--++(0,-0.2);

\draw (P1) ++ (2.7,1.1) node{$n$ half twists};

\draw (P1) ++ (0.5,0) --++(0,-0.5)to[out=-90,in=180] ++(1,-1);
\draw (P1) ++ (1,0) to[out=-90,in=180] ++(0.5,-1);

\draw (P1) ++ (1.5,-1.7) --++(2.5,0) --++(0,1.2) --++(-2.5,0) --++(0,-1.2)--cycle;

\draw (P1) ++ (4,-1.5)--++(0.3,0)  to[out=0,in=-90] ++ (0.5,0.5)--++(0,0.5);
\draw (P1) ++(4,-1) to[out=0,in=-90] ++(0.3,0.3)--++(0,0.2);

\draw (P1) ++ (2.7,-1.1) node{$n$ half twists};


\draw[preaction={draw = white, line width =5pt}] (P1) ++ (4.3,-0.5) to[out=90,in=-90] ++ (0.5,1);

\draw (P1) ++(2.5,-2) node[below]{$L(P_{1}) \cong T(2,2n+1)$};


\draw[thick] (Q) circle (2);

\draw (Q) ++(1.5,0) arc [radius =1.5, start angle = 0, end angle =180];

\draw (Q) ++ (-1.5,0) to[out=-90,in=180] ++(0.4,-0.8)to[out=0,in=-90] ++(0.2,0.4) to[out=-90,in=180] ++(0.2,-0.4) to[out=0,in=-90] ++(0.2,0.4)to[out=-90,in=180] ++(0.2,-0.4);

\draw (Q) ++ (1.5,0) to[out=-90,in=0] ++ (-0.4,-0.8) to[out=180,in=-90]++(-0.2,0.4)to[out=-90,in=0] ++ (-0.2,-0.4);

\draw (Q) ++ (0.2,-0.7) node{$\cdots$};

\draw (Q) ++(0,-1)node[below]{$n$ cusps};

\draw (Q) ++ (0,-2) node[below]{$P_2$};


\draw (Q1) ++(5,0)node[right]{$\Fix(j_{2})$};

\draw[thick] (Q1) --++(0.5,0);
\draw[thick] (Q1)++(4.5,0) --++(0.5,0);

\draw (Q1) ++(2.5,0) node{$(2n+1)$ half twists};

\draw (Q1) ++(0.5,0.5) --++(4,0) --++(0,-1) --++(-4,0)--++(0,1)--cycle;

\draw (Q1) ++ (0.5,0.3) to[out=180,in=-90] ++ (-0.3,0.5) to[out=90,in=180] ++(0.5,0.5)--++(3.6,0) to[out=0,in=90]++(0.5,-0.5)to[out=-90,in=0] ++(-0.3,-0.5);

\draw (Q1) ++ (0.5,-0.3) to[out=180,in=90] ++ (-0.3,-0.5) to[out=-90,in=180] ++(0.5,-0.5)--++(3.6,0) to[out=0,in=-90]++(0.5,0.5)to[out=90,in=0] ++(-0.3,0.5);

\draw (Q1) ++(2.5,-1.7) node[below]{$L(P_{2}) \cong T(2,2n+1)$};

\end{tikzpicture}

\caption{Two symmetries of the $(2,2n+1)$-torus knot and divides with cusps}
\label{fig:torusknot}
\end{figure}

\end{example}

\begin{figure}
\centering
\begin{tikzpicture}
\coordinate (P1) at (0,0);
\coordinate (P2) at (0,-7);

\draw[thick] (P1) --++(0,-1);
\draw[thick] (P1) ++(0,-2) --++ (0,-0.2);
\draw[thick, dashed](P1) ++ (0,-2.2) --++(0,-0.6);
\draw[thick] (P1) ++ (0,-2.8) --++(0,-0.2);
\draw[thick] (P1) ++(0,-4) --++ (0,-1);
\draw (P1) ++(0,-1) ++(0.25,0) --++(-0.5,0)--++(0,-1)--++(0.5,0)--++(0,1) --cycle;
\draw (P1) ++(0,-3) ++(0.25,0) --++(-0.5,0)--++(0,-1)--++(0.5,0)--++(0,1) --cycle;

\draw (P1) ++(0,-1.2) ++(1,0) --++(1,0)--++(0,-0.5)--++(-1,0)--++(0,0.5)--cycle;
\draw (P1) ++(0,-1.2) ++(3,0) --++(1,0)--++(0,-0.5)--++(-1,0)--++(0,0.5)--cycle;

\draw (P1) ++(0,-1.2) ++(-1,0) --++(-1,0)--++(0,-0.5)--++(1,0)--++(0,0.5)--cycle;
\draw (P1) ++(0,-1.2) ++(-3,0) --++(-1,0)--++(0,-0.5)--++(1,0)--++(0,0.5)--cycle;

\draw (P1) ++(0,-1)++(0.1,0) to[out=90,in=180] ++(0.2,0.3) to[out=0,in=180]++(0.7,-0.6);
\draw (P1) ++(0,-2)++(0.1,0) to[out=-90,in=180] ++(0.2,-0.3) to[out=0,in=180]++(0.7,0.7);

\draw (P1) ++(0,-1)++(-0.1,0) to[out=90,in=0] ++(-0.2,0.3) to[out=180,in=0]++(-0.7,-0.6);
\draw (P1) ++(0,-2)++(-0.1,0) to[out=-90,in=0] ++(-0.2,-0.3) to[out=180,in=0]++(-0.7,0.7);

\draw (P1) ++(2,0)++(0,-1.3) --++(0.2,0);
\draw[dashed] (P1) ++(2.2,-1.3) --++(0.6,0);
\draw (P1) ++(3,0)++(0,-1.3) --++(-0.2,0);

\draw (P1) ++(-2,0)++(0,-1.3) --++(-0.2,0);
\draw[dashed] (P1) ++(-2.2,-1.3) --++(-0.6,0);
\draw (P1) ++(-3,0)++(0,-1.3) --++(0.2,0);

\draw (P1) ++ (2,-1.6) to[out=0,in=30]++(0,-0.5);
\draw[dashed] (P1) ++ (2,-1.6) ++(0,-0.5) to[out=-150,in=30]++(-1.5,-0.6);
\draw (P1) ++ (2,-1.6) ++(0,-0.5) ++(-1.5,-0.6) to[out=-150,in=90] ++ (-0.4,-0.3);

\draw (P1) ++ (-2,-1.6) to[out=180,in=150]++(0,-0.5);
\draw[dashed] (P1) ++ (-2,-1.6) ++(0,-0.5) to[out=-30,in=150]++(1.5,-0.6);
\draw (P1) ++ (-2,-1.6) ++(0,-0.5) ++(1.5,-0.6) to[out=-30,in=90] ++ (0.4,-0.3);

\draw (P1) (0.1,-4) to[out=-90,in=180]++(1,-0.3) to[out=0,in=180]++(1.9,2.7);
\draw (P1) (-0.1,-4) to[out=-90,in=0]++(-1,-0.3) to[out=180,in=0]++(-1.9,2.7);

\draw (P1) ++ (-4,-1.6) to[out=180,in=120]++(0,-1) to[out=-60,in=180] ++(4,-2.2);
\draw (P1) ++ (4,-1.6) to[out=0,in=60]++(0,-1) to[out=-120,in=0] ++(-4,-2.2);

\draw (P1) ++(0,-1)++(-4,0)++(0,-0.3) to[out=180,in=180] ++(0,1) --++(8,0) to[out=0,in=0] ++(0,-1);

\draw (P1) ++(0,-1.5) node{$\alpha_1$};
\draw (P1) ++(0,-3.5) node{$\alpha_n$};
\draw (P1) ++(1.5,-1.45) node{$-c_1$};
\draw (P1) ++(-1.5,-1.45) node{$-c_1$};
\draw (P1) ++(3.5,-1.45) node{$-c_n$};
\draw (P1) ++(-3.5,-1.45) node{$-c_n$};

\draw (P1) ++(0,-5) node[below]{$I_1(\alpha_1, \ldots, \alpha_n ; c_1, \ldots, c_n)$};


\fill (P2) circle(0.06);
\draw[->] (P2) ++(-0.08,-0.2) to[out=-45,in=-60] ++(0.2,0.3);
\draw[->] (P2) ++(0.08,0.2) to[out=135,in=120] ++(-0.2,-0.3);

\draw (P2)  ++(0,0.7) --++(0.5,0);
\draw (P2) to[out=60,in=180] ++(0.5,0.4);
\draw (P2) ++ (0.5,0.3) --++(2,0) --++(0,0.5) --++(-2,0) --cycle;
\draw (P2) ++ (3.5,-0.3) --++(1,0) --++(0,-0.5) --++(-1,0) --cycle;
\draw (P2) ++ (5.1,0.3) --++(1,0) --++(0,0.5) --++(-1,0) --cycle;

\draw (P2) ++ (0,-0.7) --++(1,0);
\draw[dashed] (P2)++ (1,-0.7) --++(2,0);
\draw (P2) ++(3,-0.7) --++(0.5,0);

\draw (P2) ++ (0.5,0.3) ++ (2,0) ++ (0,0.1) to[out=0,in=120] ++(0.2,-0.2);
\draw[dashed] (P2) ++ (0.5,0.3) ++ (2,0) ++ (0,0.1) ++(0.2,-0.2) to[out=-60,in=150] ++ (0.4,-0.5);
\draw (P2) ++ (0.5,0.3) ++ (2,0) ++ (0,0.1) ++(0.2,-0.2) ++ (0.4,-0.5) to[out=-30,in=180] ++ (0.4,-0.1);

\draw (P2) ++ (2.5,0.7) --++(0.5,0);
\draw[dashed] (P2) ++(3,0.7) --++(1.5,0);
\draw (P2) ++ (4.5,0.7) --++(0.6,0);

\draw (P2) ++(4.5,-0.4) to[out=0,in=180] ++(0.6,0.8);

\draw (P2) ++(6.1,0.7) --++(0.1,0);
\draw[dashed] (P2) ++(6.2,0.7) --++(0.5,0);

\draw (P2) ++(4.5,-0.7) --++(1,0) to[out=0,in=0] ++(0.6,1.1);

\draw (P2) ++ (1.47,0.55) node{$(-1)^{n-1} a_n$};
\draw (P2) ++ (4,-0.55) node{$-a_2$};
\draw (P2) ++ (5.6,0.55) node{$a_1$};


\draw (P2)  ++(0,-0.7) --++(-0.5,0);
\draw (P2) to[out=-120,in=0] ++(-0.5,-0.4);
\draw (P2) ++ (-0.5,-0.3) --++(-2,0) --++(0,-0.5) --++(2,0) --cycle;
\draw (P2) ++ (-3.5,0.3) --++(-1,0) --++(0,0.5) --++(1,0) --cycle;
\draw (P2) ++ (-5.1,-0.3) --++(-1,0) --++(0,-0.5) --++(1,0) --cycle;

\draw (P2) ++ (0,0.7) --++(-1,0);
\draw[dashed] (P2)++ (-1,0.7) --++(-2,0);
\draw (P2) ++(-3,0.7) --++(-0.5,0);

\draw (P2) ++ (-0.5,-0.3) ++ (-2,0) ++ (0,-0.1) to[out=180,in=-60] ++(-0.2,0.2);
\draw[dashed] (P2) ++ (-0.5,-0.3) ++ (-2,0) ++ (0,-0.1) ++(-0.2,0.2) to[out=120,in=-30] ++ (-0.4,0.5);
\draw (P2) ++ (-0.5,-0.3) ++ (-2,0) ++ (0,-0.1) ++(-0.2,0.2) ++ (-0.4,0.5) to[out=150,in=0] ++ (-0.4,0.1);

\draw (P2) ++ (-2.5,-0.7) --++(-0.5,0);
\draw[dashed] (P2) ++(-3,-0.7) --++(-1.5,0);
\draw (P2) ++ (-4.5,-0.7) --++(-0.6,0);

\draw (P2) ++(-4.5,0.4) to[out=180,in=0] ++(-0.6,-0.8);

\draw (P2) ++(-6.1,-0.7) --++(-0.1,0);
\draw[dashed] (P2) ++(-6.2,-0.7) --++(-0.5,0);

\draw (P2) ++(-4.5,0.7) --++(-1,0) to[out=180,in=180] ++(-0.6,-1.1);

\draw (P2) ++ (-1.47,-0.55) node{$(-1)^{n-1} a_n$};
\draw (P2) ++ (-4,0.55) node{$-a_2$};
\draw (P2) ++ (-5.6,-0.55) node{$a_1$};

\draw (P2) ++(0,-1) node[below]{$I_2 (a_1, \ldots, a_n)$};

\draw (P2) ++(0,-2.25) node{$=$};
\draw (P2) ++ (-0.3,-2) --++(-1,0) --++(0,-0.5) --++(1,0) --cycle;
\draw (P2) ++(-0.8,-2.25) node{$b$};

\draw (P2) ++(0.3,-2.1) --++(0.3,0);
\draw (P2) ++(0.3,-2.4) --++(0.3,0);

\draw (P2) ++ (0.6,-2.4) to[out=0,in=180] ++(0.3,0.3);
\draw[preaction={draw=white, line width = 3pt}] (P2) ++ (0.6,-2.1) to[out=0,in=180] ++(0.3,-0.3);

\draw (P2) ++ (0.9,-2.4) to[out=0,in=180] ++(0.3,0.3);
\draw[preaction={draw=white, line width = 3pt}] (P2) ++ (0.9,-2.1) to[out=0,in=180] ++(0.3,-0.3);

\draw (P2) ++(1.5,-2.25) node{$\cdots$};
\draw (P2) ++ (1.8,-2.4) to[out=0,in=180] ++(0.3,0.3) --++(0.3,0);
\draw[preaction={draw=white, line width = 3pt}] (P2) ++ (1.8,-2.1) to[out=0,in=180] ++(0.3,-0.3) --++(0.3,0);

\draw (P2) ++(1.4,-2.6) node{{\footnotesize $b$ right half twists}};

\end{tikzpicture}
\caption{Knot diagrams of $I_1 (\alpha_1, \ldots, \alpha_n ; c_1 ,\ldots, c_n)$ and $I_2 (a_1 ,\ldots, a_n)$}
\label{fig:sakuma}
\end{figure}
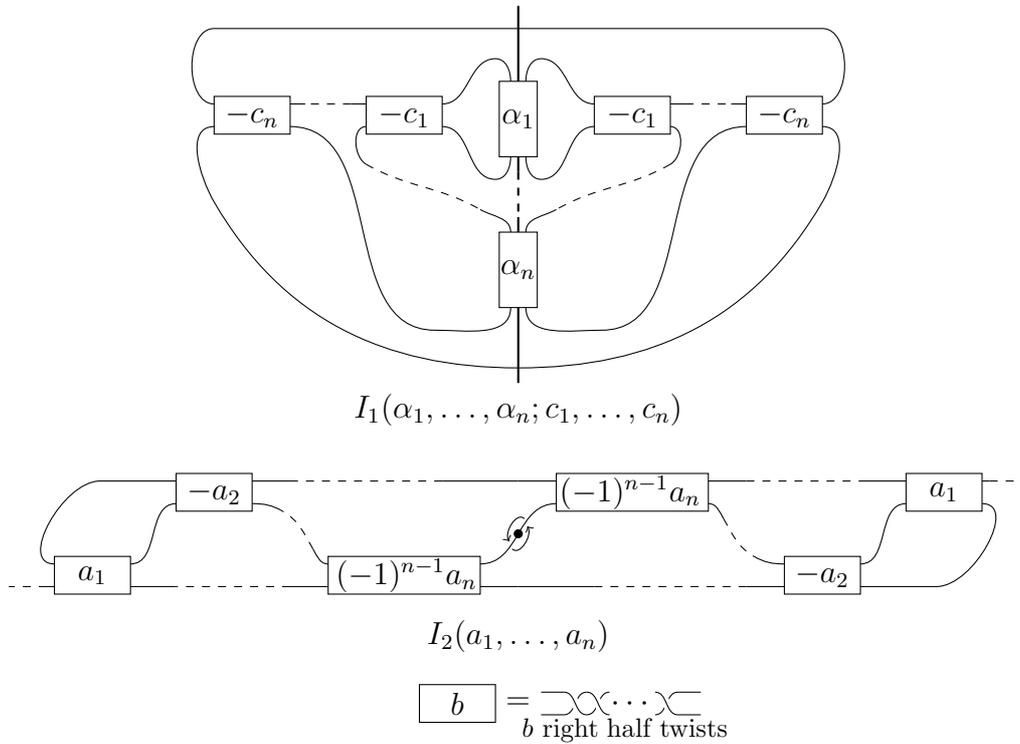

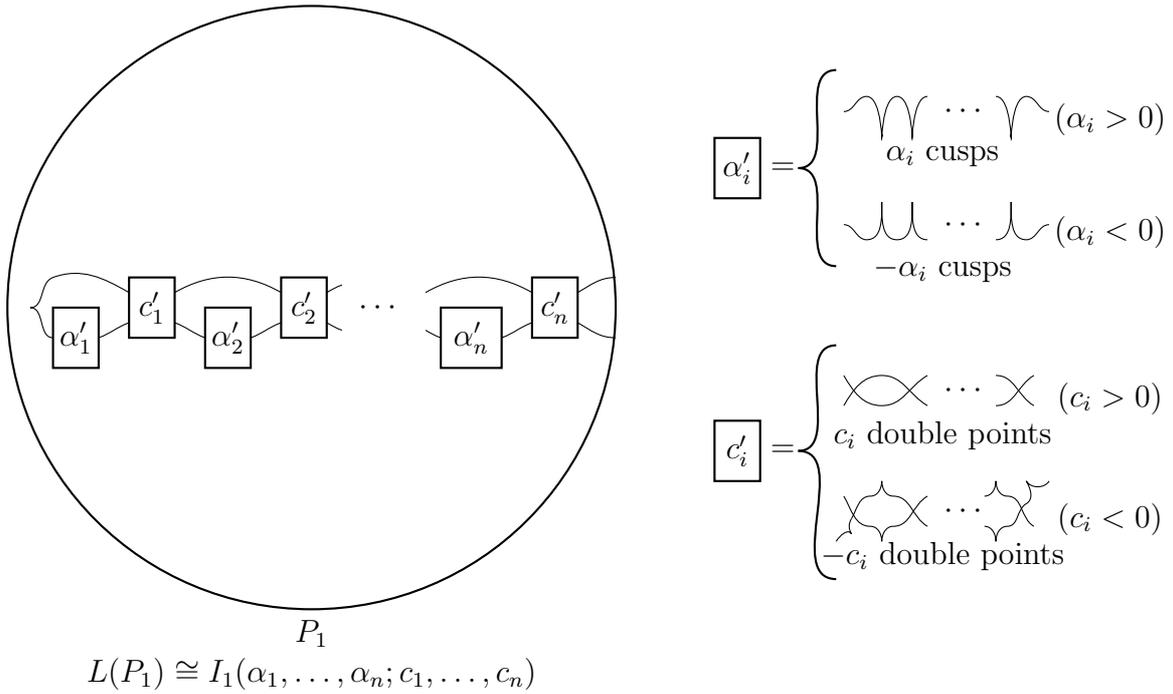
\begin{figure}
\centering
\begin{tikzpicture}
\coordinate (P1) at (0,0);
\coordinate (P2) at (-3.4,0);
\coordinate (Q) at (5.5,1.5);
\coordinate (R) at (5.5,-1.5);

\draw[thick] (P1) circle (4);
\draw (P1) ++(0,-4) node[below]{$P_1$};
\draw (P1) ++(0,-4.5) node[below] {$L(P_1) \cong I_1 (\alpha_1, \ldots, \alpha_n ; c_1, \ldots, c_n)$};
 
\draw[thick] (P2) --++(0.6,0) --++(0,-0.8) --++(-0.6,0) --++(0,0.8) --cycle;
\draw[thick] (P2) ++ (1,0.4) --++(0.6,0) --++(0,-0.8) --++(-0.6,0) --++(0,0.8) --cycle;
\draw[thick] (P2) ++ (2,-0) --++(0.6,0) --++(0,-0.8) --++(-0.6,0) --++(0,0.8) --cycle;
\draw[thick] (P2) ++ (3,0.4) --++(0.6,0) --++(0,-0.8) --++(-0.6,0) --++(0,0.8) --cycle;
\draw[thick] (P2) ++ (5.1,0) --++(0.8,0) --++(0,-0.8) --++(-0.8,0) --++(0,0.8) --cycle;
\draw[thick] (P2) ++ (6.3,0.4) --++(0.6,0) --++(0,-0.8) --++(-0.6,0) --++(0,0.8) --cycle;

\draw (P2) ++(0.3,-0.4) node{$\alpha'_{1}$};
\draw (P2) ++(1,0.4)++(0.3,-0.4) node{$c'_{1}$};
\draw (P2) ++ (2,-0)++(0.3,-0.4) node{$\alpha'_{2}$};
\draw (P2) ++(3,0.4)++(0.3,-0.4) node{$c'_{2}$};
\draw (P2) ++ (5.2,-0)++(0.3,-0.4) node{$\alpha'_{n}$};
\draw (P2) ++(6.3,0.4)++(0.3,-0.4) node{$c'_{n}$};

\draw (P2) ++ (0,-0.4) to [out=180,in=0] ++ (-0.3,0.4) to [out=0,in=200] ++ (0.3,0.4) to[out=20,in=150] ++ (1,-0.2);

\draw (P2) ++ (1.6,0.2) to[out=30,in=150] ++(1.4,0);
\draw (P2) ++ (3.6,0.2) to[out=30, in=200] ++ (0.2,0.1);
\draw (P2) ++ (4.9,0.2) to[out=30,in=150] ++(1.4,0);
\draw (P2) ++ (6.9,0.2) to[out=30,in=180] ++(0.5,0.2);

\draw (P2) ++ (0.6,-0.4) to[out=20,in=200]++ (0.4,0.2);
\draw (P2) ++ (1.6,-0.2) to[out=-20,in=160]++ (0.4,-0.2);
\draw (P2) ++ (2.6,-0.4) to[out=20,in=200]++ (0.4,0.2);
\draw (P2) ++ (3.6,-0.2) to[out=-30,in=160]++ (0.2,-0.1);
\draw (P2) ++ (4.9,-0.3) to[out=-30,in=160]++ (0.2,-0.1);
\draw (P2) ++ (5.9,-0.4) to[out=20,in=200]++ (0.4,0.2);
\draw (P2) ++ (6.9,-0.2) to[out=-30,in=180] ++(0.5,-0.2);

\draw (P2) ++ (4.3,0) node{$\cdots$};


\draw[thick] (Q) ++(-0.2,0.75) --++(0.6,0) --++(0,-0.8) --++(-0.6,0) --++(0,0.8) --cycle;
\draw (Q) ++(-0.2,0.75)++ (0.9,-0.4) node{$=$};
\draw (Q) ++(-0.2,0.75)++ (0.3,-0.4) node{$\alpha'_{i}$};
\draw[thick] (Q) ++ (-0.2,0.75) ++ (1.1,-0.4) to[out=0,in=180] ++(0.5,1.3);
\draw[thick] (Q) ++ (-0.2,0.75) ++ (1.1,-0.4) to[out=0,in=180] ++(0.5,-1.3);

\draw (Q) ++(0,1.5) ++(1.5,-0.4) to[out=0,in=180] ++ (0.3,0.2) to [out=0,in=-90] ++(0.2,-0.5) to [out=-90,in=180] ++(0.2,0.5) to [out=0,in=-90] ++(0.2,-0.5) to [out=-90,in=180] ++(0.2,0.5) ; 
\draw (Q) ++(0,1.5) ++ (3.5,-0.2) to [out=0,in=-90] ++(0.2,-0.5) to [out=-90,in=180] ++(0.2,0.5) to[out=0,in=180] ++(0.3,-0.2);

\draw (Q) ++(0,1.5) ++ (3.1,-0.4) node{$\cdots$} ;

\draw (Q) ++(0,1.5)++ (2.8,-1) node{$\alpha_i$ cusps};
\draw (Q) ++(0,1.5)++ (5,-0.50) node{$(\alpha_i >0)$};

\draw (Q) ++(1.5,-0.4) to[out=0,in=180] ++ (0.3,-0.2) to [out=0,in=-90] ++(0.2,0.5) to [out=-90,in=180] ++(0.2,-0.5) to [out=0,in=-90] ++(0.2,0.5) to [out=-90,in=180] ++(0.2,-0.5) ; 
\draw (Q) ++ (3.5,-0.6) to [out=0,in=-90] ++(0.2,0.5) to [out=-90,in=180] ++(0.2,-0.5) to[out=0,in=180] ++(0.3,0.2);

\draw (Q) ++ (3.1,-0.4) node{$\cdots$} ;

\draw (Q) ++ (2.8,-1) node{$-\alpha_i$ cusps};
\draw (Q) ++ (5,-0.50) node{$(\alpha_i <0)$};

\draw[thick] (R) ++(-0.2,0)--++(0.6,0) --++(0,-0.8) --++(-0.6,0) --++(0,0.8) --cycle;
\draw (R)++(-0.2,0) ++ (0.9,-0.4) node{$=$};
\draw (R)++(-0.2,0) ++ (0.3,-0.4) node{$c'_{i}$};
\draw[thick] (R)++(-0.2,0)  ++ (1.1,-0.4) to[out=0,in=180] ++(0.5,1.4);
\draw[thick] (R) ++(-0.2,0) ++ (1.1,-0.4) to[out=0,in=180] ++(0.5,-1.7);

\draw (R) ++(0,0.8)++(1.5,-0.6) to [out= 60, in=180] ++(0.5,0.4) to[out=0,in=150] ++(0.6,-0.4); 
\draw (R)++(0,0.8) ++ (1.5,-0.2)  to [out= -60, in=180] ++(0.5,-0.4) to[out=0,in=-150] ++(0.6,0.4); 
\draw (R)++(0,0.8) ++ (3.5,-0.6) to [out= 0, in=210] ++(0.5,0.4);
\draw (R) ++(0,0.8)++ (3.5,-0.2) to [out= 0, in=150] ++(0.5,-0.4);

\draw (R) ++(0,0.8)++ (3.1,-0.4) node{$\cdots$} ;

\draw (R) ++(0,0.8)++ (2.8,-1) node{$c_i$ double points};
\draw (R)++(0,0.8) ++ (5,-0.50) node{$(c_i >0)$};

\draw (R) ++ (0,-0.8) ++ (1.4, -0.8)to[out=60,in=120] ++(0.2,0.1)to[out=120,in=-120] ++(0,0.2);
\draw (R) ++(0,-0.8)++(1.6,-0.5) to [out= 60, in=180] ++(0.25,0.3)to[out=0,in=-90]++(0.15,0.2)to[out=-90,in=180]++(0.15,-0.2) to[out=0,in=150] ++(0.45,-0.4);

\draw (R)++(0,-0.8) ++ (1.5,-0.2)  to [out= -60, in=180] ++(0.35,-0.4)to[out=0,in=90] ++(0.15,-0.2)to[out=90,in=180] ++ (0.15,0.2) to[out=0,in=-150] ++(0.45,0.4);
 
\draw (R)++(0,-0.8) ++ (3.35,-0.6)to[out=0,in=90] ++(0.15,-0.2)to[out=90,in=180] ++ (0.15,0.2)  to [out= 0, in=210] ++(0.25,0.3) to[out=30,in=-60]++(0,0.3) to[out=-60, in=-120] ++(0.3,0);

\draw (R) ++(0,-0.8)++ (3.35,-0.2)to[out=0,in=-90]++(0.15,0.2)to[out=-90,in=180]++(0.15,-0.2) to [out= 0, in=150] ++(0.35,-0.4);

\draw (R) ++(0,-0.8)++ (3.1,-0.4) node{$\cdots$} ;

\draw (R) ++(0,-0.8)++ (2.8,-1) node{$-c_i$ double points};
\draw (R)++(0,-0.8) ++ (5,-0.50) node{$(c_i <0)$};

\end{tikzpicture}
\caption{A divide whose associated link is $I_1 (\alpha_1, \ldots, \alpha_n ; c_1 ,\ldots, c_n)$ }
\label{fig:sakumaI1}
\end{figure}

\begin{figure}[htbp]
\centering

\begin{tikzpicture}
\coordinate (P1) at (0,0);
\coordinate (P2) at (-3.4,0);
\coordinate (Q) at (5.5,1.5);
\coordinate (R) at (5.5,-1.5);

\draw[thick] (P1) circle (4);
\draw[thick] (P2) --++(0.6,0) --++(0,-0.8) --++(-0.6,0) --++(0,0.8) --cycle;
\draw[thick] (P2) ++ (1,0.4) --++(0.6,0) --++(0,-0.8) --++(-0.6,0) --++(0,0.8) --cycle;
\draw[thick] (P2) ++ (2,-0) --++(0.6,0) --++(0,-0.8) --++(-0.6,0) --++(0,0.8) --cycle;
\draw[thick] (P2) ++ (3,0.4) --++(0.6,0) --++(0,-0.8) --++(-0.6,0) --++(0,0.8) --cycle;
\draw[thick] (P2) ++ (5.1,0) --++(0.8,0) --++(0,-0.8) --++(-0.8,0) --++(0,0.8) --cycle;
\draw[thick] (P2) ++ (6.3,0.4) --++(0.6,0) --++(0,-0.8) --++(-0.6,0) --++(0,0.8) --cycle;

\draw (P2) ++(0.3,-0.4) node{$b_{1}$};
\draw (P2) ++(1,0.4)++(0.3,-0.4) node{$c_{2}$};
\draw (P2) ++ (2,-0)++(0.3,-0.4) node{$b_{3}$};
\draw (P2) ++(3,0.4)++(0.3,-0.4) node{$c_{4}$};
\draw (P2) ++ (5.2,-0)++(0.3,-0.4) node{$b_{n-1}$};
\draw (P2) ++(6.3,0.4)++(0.3,-0.4) node{$c_{n}$};

\draw (P2) ++ (0,-0.4) to [out=180,in=0] ++ (-0.3,0.4) to [out=0,in=200] ++ (0.3,0.4) to[out=20,in=150] ++ (1,-0.2);

\draw (P2) ++ (1.6,0.2) to[out=30,in=150] ++(1.4,0);
\draw (P2) ++ (3.6,0.2) to[out=30, in=200] ++ (0.2,0.1);
\draw (P2) ++ (4.9,0.2) to[out=30,in=150] ++(1.4,0);
\draw (P2) ++ (6.9,0.2) to[out=30,in=180] ++(0.5,0.2);

\draw (P2) ++ (0.6,-0.4) to[out=20,in=200]++ (0.4,0.2);
\draw (P2) ++ (1.6,-0.2) to[out=-20,in=160]++ (0.4,-0.2);
\draw (P2) ++ (2.6,-0.4) to[out=20,in=200]++ (0.4,0.2);
\draw (P2) ++ (3.6,-0.2) to[out=-30,in=160]++ (0.2,-0.1);
\draw (P2) ++ (4.9,-0.3) to[out=-30,in=160]++ (0.2,-0.1);
\draw (P2) ++ (5.9,-0.4) to[out=20,in=200]++ (0.4,0.2);
\draw (P2) ++ (6.9,-0.2) to[out=-30,in=180] ++(0.5,-0.2);

\draw (P2) ++ (4.3,0) node{$\cdots$};

\draw[thick] (Q) --++(0.6,0) --++(0,-0.8) --++(-0.6,0) --++(0,0.8) --cycle;
\draw (Q) ++ (1,-0.4) node{$=$};

\draw (Q) ++(1.5,-0.4) to[out=-60,in=180] ++ (0.3,-0.2) to [out=0,in=-90] ++(0.2,0.5) to [out=-90,in=180] ++(0.2,-0.5) to [out=0,in=-90] ++(0.2,0.5) to [out=-90,in=180] ++(0.2,-0.5) ; 
\draw (Q) ++ (3.5,-0.6) to [out=0,in=-90] ++(0.2,0.5) to [out=-90,in=180] ++(0.2,-0.5) to[out=0,in=240] ++(0.3,0.2);

\draw (Q) ++ (3.1,-0.4) node{$\cdots$} ;

\draw (Q) ++ (2.8,-1) node{$a_i$ cusps};
\draw (Q) ++ (0.3,-0.4) node{$b_{i}$};

\draw[thick] (R) --++(0.6,0) --++(0,-0.8) --++(-0.6,0) --++(0,0.8) --cycle;
\draw (R) ++ (1,-0.4) node{$=$};

\draw (R) ++(1.5,-0.6) to [out= 60, in=180] ++(0.5,0.4) to[out=0,in=150] ++(0.6,-0.4); 
\draw (R) ++ (1.5,-0.2)  to [out= -60, in=180] ++(0.5,-0.4) to[out=0,in=-150] ++(0.6,0.4); 
\draw (R) ++ (3.5,-0.6) to [out= 0, in=210] ++(0.5,0.4);
\draw (R) ++ (3.5,-0.2) to [out= 0, in=150] ++(0.5,-0.4);

\draw (R) ++ (3.1,-0.4) node{$\cdots$} ;

\draw (R) ++ (2.8,-1) node{$a_i /2$ double points};
\draw (R) ++ (0.3,-0.4) node{$c_{i}$};

\end{tikzpicture}

\caption{A divide with cusps whose link is $C(a_{1} , \cdots , a_{n})$}
\label{fig:2bridge}
\end{figure}
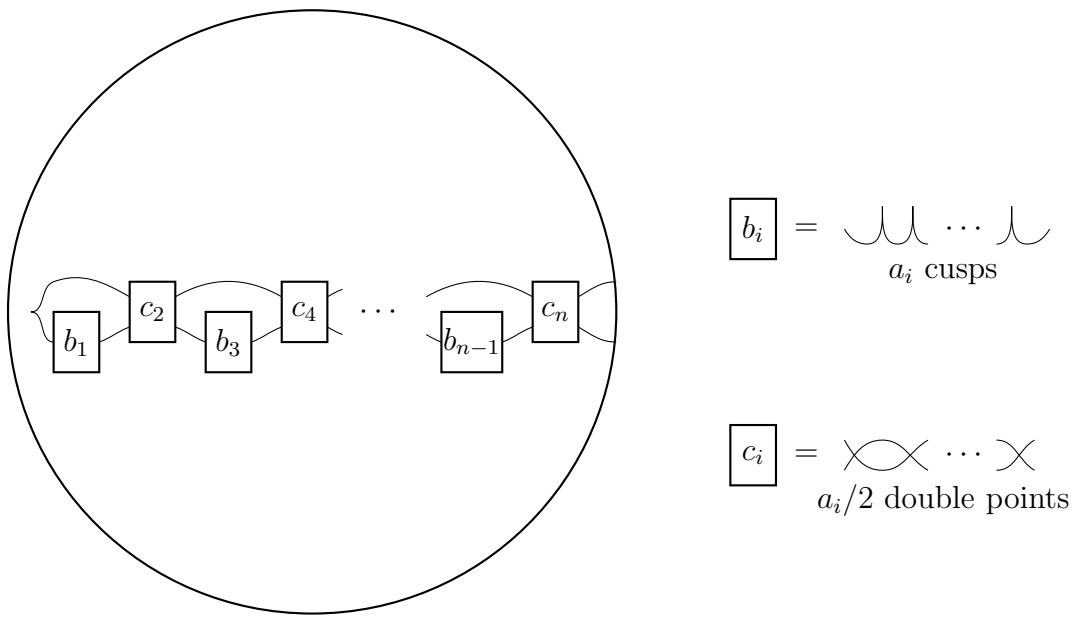

\begin{example}
In \cite{sakuma}, Sakuma gave some explicit descriptions of strongly invertible knots, and they are written as $I_{1} (\alpha_1, \ldots, \alpha_{n}; c_{1} , \ldots, c_{n} )$ and $I_{2} (a_{1} , \ldots, a_{n})$ (see Figure \ref{fig:sakuma}), where each $\alpha_i, c_i, a_i$ is a non-negative integer and $a_i$ is even. By using this diagram, we can see that a divide with cusps whose link is such a strongly invertible knot is described in Figure \ref{fig:sakumaI1} and \ref{fig:sakumaI2}. 

A $2$-bridge knot can be described by using this notation and it is known that a $2$-bridge knot has exactly two symmetries such that it is strongly invertible (Proposition 3.6 in \cite{sakuma}). A $2$-bridge knot is expressed as by using a Conway notation $C(a_{1} , \ldots, a_{n})$ for even numbers $a_{i}$. One of the strongly invertible knots defined by $C(a_{1}, \ldots, a_{n})$ is described in Figure \ref{fig:conway}, and the corresponding divide with cusps whose link is $C(a_{1}, \ldots, a_{n})$ is as in Figure \ref{fig:2bridge}.
\end{example}

\begin{figure}
\centering
\begin{tikzpicture}


\coordinate (A1) at (0,0);
\coordinate (A2) at (-3.2,0);
\coordinate (B) at (5.5,-1.5);
\coordinate (C) at (2.5,-6.5);

\coordinate (D1) at (9,0);
\coordinate (D2) at (5.8,0);
\coordinate (E) at (14.5,-1.5);
\coordinate (F) at (12,-6.5);

\draw[thick] (A1) circle (4);
\draw (A1) ++(0,-4) node[below]{$P_2$};
\draw (A1) ++(0,-4.5) node[below] {$L(P_2) \cong I_2 (a_1, \ldots, a_n)$ ($n$ is odd)};
 
\draw[thick] (A2) ++ (0,-0.1) --++(0.6,0) --++(0,0.8) --++(-0.6,0) --++(0,-0.8) --cycle;
\draw[thick] (A2) ++ (0.9,0.1) --++(0.8,0) --++(0,-0.8) --++(-0.8,0) --++(0,0.8) --cycle;
\draw[thick] (A2) ++ (2,-0.1) --++(0.6,0) --++(0,0.8) --++(-0.6,0) --++(0,-0.8) --cycle;
\draw[thick] (A2) ++ (3.6,0.1) --++(1.1,0) --++(0,-0.8) --++(-1.1,0) --++(0,0.8) --cycle;
\draw[thick] (A2) ++ (4.9,-0.1) --++(1.5,0) --++(0,0.8) --++(-1.5,0) --++(0,-0.8) --cycle;

\draw (A2) ++ (0.3,0.3) node{$a_1$};
\draw (A2) ++ (1.3,-0.3) node{$-a_2$};
\draw (A2) ++ (2.3,0.3) node{$a_3$};
\draw (A2) ++ (4.15,-0.3) node{\small{$-a_{n-1}$}};
\draw (A2) ++ (5.55,0.3) node{\small{$a_{n}+1$}};

\draw (A2) ++(0,0.6) --++(-0.77,0);
\draw (A2) to[out=180,in=180]++(0,-0.6) --++(0.2,0) to[out=0,in=90] ++(0.1,-0.2) to[out=90,in=180] ++(0.1,0.2) --++(0.2,0); 

\draw (A2) ++ (0.6,0.6) --++(0.3,0);
\draw (A2) ++ (0.6,0) --++(0.3,0);
\draw (A2) ++ (0.6,-0.6) --++(0.3,0);

\draw (A2) ++ (0.9,0.6) --++(0.8,0);

\draw (A2) ++ (1.7,0.6) --++(0.3,0);
\draw (A2) ++ (1.7,0) --++(0.3,0);
\draw (A2) ++ (1.7,-0.6) --++(0.3,0);

\draw (A2) ++ (2.0,-0.6) --++(0.6,0);

\draw (A2) ++ (2.6,0.6) --++(0.1,0);
\draw (A2) ++ (2.6,0) --++(0.1,0);
\draw (A2) ++ (2.6,-0.6) --++(0.1,0);

\draw (A2) ++ (3.2,0) node{$\cdots$};

\draw (A2) ++ (3.5,0.6) --++(0.1,0);
\draw (A2) ++ (3.5,0) --++(0.1,0);
\draw (A2) ++ (3.5,-0.6) --++(0.1,0);

\draw (A2) ++ (3.6,0.6) --++(1.1,0);

\draw (A2) ++ (4.7,0.6) --++(0.2,0);
\draw (A2) ++ (4.7,0) --++(0.2,0);
\draw (A2) ++ (4.7,-0.6) --++(0.2,0);

\draw (A2) ++ (4.9,-0.6) --++(1.5,0);

\draw (A2) ++(6.4,0.6) --++(0.75,0); 
\draw (A2) ++ (6.4,0) to[out=0,in=0] ++(0,-0.6);

\draw[thick] (D1) circle (4);
\draw (D1) ++(0,-4) node[below]{$P_2$};
\draw (D1) ++(0,-4.5) node[below] {$L(P_2) \cong I_2 (a_1, \ldots, a_n) $ ($n$ is even)};
 
\draw[thick] (D2) ++ (0,-0.1) --++(0.6,0) --++(0,0.8) --++(-0.6,0) --++(0,-0.8) --cycle;
\draw[thick] (D2) ++ (0.9,0.1) --++(0.8,0) --++(0,-0.8) --++(-0.8,0) --++(0,0.8) --cycle;
\draw[thick] (D2) ++ (2,-0.1) --++(0.6,0) --++(0,0.8) --++(-0.6,0) --++(0,-0.8) --cycle;
\draw[thick] (D2) ++ (3.7,-0.1) --++(0.9,0) --++(0,0.8) --++(-0.9,0) --++(0,-0.8) --cycle;
\draw[thick] (D2) ++ (4.9,0.1) --++(1.7,0) --++(0,-0.8) --++(-1.7,0) --++(0,0.8) --cycle;

\draw (D2) ++ (0.3,0.3) node{$a_1$};
\draw (D2) ++ (1.3,-0.3) node{$-a_2$};
\draw (D2) ++ (2.3,0.3) node{$a_3$};
\draw (D2) ++ (4.15,0.3) node{\small{$a_{n-1}$}};
\draw (D2) ++ (5.75,-0.3) node{\small{$-(a_{n}+1)$}};

\draw (D2) ++(0,0.6) --++(-0.77,0);
\draw (D2) to[out=180,in=180]++(0,-0.6) --++(0.2,0) to[out=0,in=90] ++(0.1,-0.2) to[out=90,in=180] ++(0.1,0.2) --++(0.2,0); 

\draw (D2) ++ (0.6,0.6) --++(0.3,0);
\draw (D2) ++ (0.6,0) --++(0.3,0);
\draw (D2) ++ (0.6,-0.6) --++(0.3,0);

\draw (D2) ++ (0.9,0.6) --++(0.8,0);

\draw (D2) ++ (1.7,0.6) --++(0.3,0);
\draw (D2) ++ (1.7,0) --++(0.3,0);
\draw (D2) ++ (1.7,-0.6) --++(0.3,0);

\draw (D2) ++ (2.0,-0.6) --++(0.6,0);

\draw (D2) ++ (2.6,0.6) --++(0.1,0);
\draw (D2) ++ (2.6,0) --++(0.1,0);
\draw (D2) ++ (2.6,-0.6) --++(0.1,0);

\draw (D2) ++ (3.1,0) node{$\cdots$};

\draw (D2) ++ (3.6,0.6) --++(0.1,0);
\draw (D2) ++ (3.6,0) --++(0.1,0);
\draw (D2) ++ (3.6,-0.6) --++(0.1,0);

\draw (D2) ++ (3.7,-0.6) --++(0.9,0);

\draw (D2) ++ (4.6,0.6) --++(0.3,0);
\draw (D2) ++ (4.6,0) --++(0.3,0);
\draw (D2) ++ (4.6,-0.6) --++(0.3,0);

\draw (D2) ++ (4.9,0.6) --++(0.4,0) to[out=0,in=-90] ++(0.1,0.2)to[out=-90.,in=180] ++(0.1,-0.2) --++(0.5,0)to[out=0,in=-90] ++(0.1,0.2)to[out=-90.,in=180] ++(0.1,-0.2) --++(0.4,0);

\draw (D2) ++(6.6,-0.6) --++(0.55,0); 
\draw (D2) ++ (6.6,0) to[out=0,in=0] ++(0,0.6);

\draw[thick] (C) ++(-0.2,0)--++(0.6,0) --++(0,-0.8) --++(-0.6,0) --++(0,0.8) --cycle;
\draw (C)++(-0.2,0) ++ (0.9,-0.4) node{$=$};
\draw (C)++(-0.2,0) ++ (0.3,-0.4) node{$a$};
\draw[thick] (C)++(-0.2,0)  ++ (1.1,-0.4) to[out=0,in=180] ++(0.5,1.4);
\draw[thick] (C) ++(-0.2,0) ++ (1.1,-0.4) to[out=0,in=180] ++(0.5,-1.7);

\draw (C) ++(0,0.8)++(1.5,-0.6) to [out= 60, in=180] ++(0.5,0.4) to[out=0,in=150] ++(0.6,-0.4); 
\draw (C)++(0,0.8) ++ (1.5,-0.2)  to [out= -60, in=180] ++(0.5,-0.4) to[out=0,in=-150] ++(0.6,0.4); 
\draw (C)++(0,0.8) ++ (3.5,-0.6) to [out= 0, in=210] ++(0.5,0.4);
\draw (C) ++(0,0.8)++ (3.5,-0.2) to [out= 0, in=150] ++(0.5,-0.4);

\draw (C) ++(0,0.8)++ (3.1,-0.4) node{$\cdots$} ;

\draw (C) ++(0,0.8)++ (2.8,-1) node{$a$ double points};
\draw (C)++(0,0.8) ++ (5,-0.50) node{$(a >0)$};

\draw (C) ++ (0,-0.8) ++ (1.4, -0.8)to[out=60,in=120] ++(0.2,0.1)to[out=120,in=-120] ++(0,0.2);
\draw (C) ++(0,-0.8)++(1.6,-0.5) to [out= 60, in=180] ++(0.25,0.3)to[out=0,in=-90]++(0.15,0.2)to[out=-90,in=180]++(0.15,-0.2) to[out=0,in=150] ++(0.45,-0.4);

\draw (C)++(0,-0.8) ++ (1.5,-0.2)  to [out= -60, in=180] ++(0.35,-0.4)to[out=0,in=90] ++(0.15,-0.2)to[out=90,in=180] ++ (0.15,0.2) to[out=0,in=-150] ++(0.45,0.4);
 
\draw (C)++(0,-0.8) ++ (3.35,-0.6)to[out=0,in=90] ++(0.15,-0.2)to[out=90,in=180] ++ (0.15,0.2)  to [out= 0, in=210] ++(0.25,0.3) to[out=30,in=-60]++(0,0.3) to[out=-60, in=-120] ++(0.3,0);

\draw (C) ++(0,-0.8)++ (3.35,-0.2)to[out=0,in=-90]++(0.15,0.2)to[out=-90,in=180]++(0.15,-0.2) to [out= 0, in=150] ++(0.35,-0.4);

\draw (C) ++(0,-0.8)++ (3.1,-0.4) node{$\cdots$} ;

\draw (C) ++(0,-0.8)++ (2.8,-1) node{$-a$ double points};
\draw (C)++(0,-0.8) ++ (5,-0.50) node{$(a <0)$};

\end{tikzpicture}
\caption{A divide whose associated link is $I_2 (a_1, \ldots, a_n)$}
\label{fig:sakumaI2}
\end{figure}


\begin{figure}[htbp]
\centering

\begin{tikzpicture}
\coordinate (P1) at (0,0);
\coordinate (P2) at (0,-3);
\draw (P1) to[out=180,in=-90] ++ (-0.4,0.3) to[out=90,in=180] ++(0.4,0.3);
\draw (P1) ++ (0,-1.2) to[out=180,in=-90] ++ (-0.4,0.3) to[out=90,in=180] ++(0.4,0.3);

\draw[thick] (P1) ++ (0,0.1) --++(0,-0.8) --++(1,0) --++(0,0.8) --++(-1,0) --cycle;
\draw[thick] (P1) ++ (0,0.6) ++ (0,0.1) ++ (1.3,0) --++(0,-0.8) --++(1,0) --++(0,0.8) --++(-1,0) --cycle;
\draw[thick] (P1) ++ (2.6,0)++ (0,0.1) --++(0,-0.8) --++(1,0) --++(0,0.8) --++(-1,0) --cycle;
\draw[thick] (P1) ++ (5.3,0)++ (0,0.1) --++(0,-0.8) --++(1.2,0) --++(0,0.8) --++(-1.2,0) --cycle;
\draw[thick] (P1) ++ (0,0.6) ++ (0,0.1) ++ (6.8,0) --++(0,-0.8) --++(1,0) --++(0,0.8) --++(-1,0) --cycle;

\draw (P1) ++ (0,0.6) --++(1.3,0);
\draw (P1) ++ (2.3,0.6) --++(1.7,0);
\draw (P1) ++ (5.2,0.6) --++(1.6,0);
\draw (P1) ++ (1,0) --++(0.3,0);
\draw (P1) ++ (2.3,0) --++(0.3,0);
\draw (P1) ++ (3.6,0) --++(0.4,0);
\draw (P1) ++ (5.2,0) --++(0.1,0);
\draw (P1) ++ (6.5,0) --++(0.3,0);
\draw (P1) ++ (1,-0.6) --++(1.6,0);
\draw (P1) ++ (3.6,-0.6) --++(0.4,0);
\draw (P1) ++ (5.2,-0.6) --++(0.1,0);
\draw (P1) ++ (6.5,-0.6) --++(1.3,0);
\draw (P1) ++ (0,-1.2) --++(4,0);
\draw (P1) ++ (5.2,-1.2) --++(2.6,0);

\draw (P1) ++ (7.8,0) to[out=0,in=90] ++ (0.4,-0.3) to[out=-90,in=0] ++(-0.4,-0.3);
\draw (P1) ++ (7.8,0) ++ (0,-1.2) to[out=0,in=-90] ++ (0.9,0.9) to[out=90,in=0] ++(-0.9,0.9);

\draw (P1) ++(0.5,-0.3) node{$-a_1$};
\draw (P1) ++(3.1,-0.3) node{$-a_3$};
\draw (P1) ++(5.9,-0.3) node{$-a_{n-1}$};

\draw (P1) ++ (1.8,0.3) node{$a_2$};
\draw (P1) ++ (7.3,0.3) node{$a_n$};

\draw (P1) ++ (4.6,0.3) node{$\cdots$};
\draw (P1) ++ (4.6,-0.9) node{$\cdots$};

\draw[->, very thick] (P1) ++ (4,-1.6) --++(0,-0.5);


\draw (P2) to[out=180,in=-90] ++ (-0.4,0.3) to[out=90,in=180] ++(0.4,0.3);
\draw (P2) ++ (0,-1.2) to[out=180,in=-90] ++ (-0.4,0.3) to[out=90,in=180] ++(0.4,0.3);

\draw[thick] (P2) ++ (0,0.1) --++(0,-0.8) --++(1,0) --++(0,0.8) --++(-1,0) --cycle;
\draw[thick] (P2) ++ (0,0.6) ++ (0,0.1) ++ (1.3,0) --++(0,-0.8) --++(1,0) --++(0,0.8) --++(-1,0) --cycle;
\draw[thick] (P2) ++ (0,-0.6) ++ (0,0.1) ++ (1.3,0) --++(0,-0.8) --++(1,0) --++(0,0.8) --++(-1,0) --cycle;
\draw[thick] (P2) ++ (2.6,0)++ (0,0.1) --++(0,-0.8) --++(1,0) --++(0,0.8) --++(-1,0) --cycle;
\draw[thick] (P2) ++ (5.3,0)++ (0,0.1) --++(0,-0.8) --++(1.2,0) --++(0,0.8) --++(-1.2,0) --cycle;
\draw[thick] (P2) ++ (0,0.6) ++ (0,0.1) ++ (6.8,0) --++(0,-0.8) --++(1,0) --++(0,0.8) --++(-1,0) --cycle;
\draw[thick] (P2) ++ (0,-0.6) ++ (0,0.1) ++ (6.8,0) --++(0,-0.8) --++(1,0) --++(0,0.8) --++(-1,0) --cycle;

\draw (P2) ++ (0,0.6) --++(1.3,0);
\draw (P2) ++ (2.3,0.6) --++(1.7,0);
\draw (P2) ++ (5.2,0.6) --++(1.6,0);
\draw (P2) ++ (1,0) --++(0.3,0);
\draw (P2) ++ (2.3,0) --++(0.3,0);
\draw (P2) ++ (3.6,0) --++(0.4,0);
\draw (P2) ++ (5.2,0) --++(0.1,0);
\draw (P2) ++ (6.5,0) --++(0.3,0);
\draw (P2) ++ (1,-0.6) --++(0.3,0);
\draw (P2) ++ (2.3,-0.6) --++(0.3,0);
\draw (P2) ++ (3.6,-0.6) --++(0.4,0);
\draw (P2) ++ (5.2,-0.6) --++(0.1,0);
\draw (P2) ++ (6.5,-0.6) --++(0.3,0);

\draw (P2) ++ (0,-1.2) --++(1.3,0);
\draw (P2) ++ (2.3,-1.2) --++(1.7,0);
\draw (P2) ++ (5.2,-1.2) --++(1.6,0);

\draw (P2) ++ (7.8,0) to[out=0,in=90] ++ (0.4,-0.3) to[out=-90,in=0] ++(-0.4,-0.3);
\draw (P2) ++ (7.8,0) ++ (0,-1.2) to[out=0,in=-90] ++ (0.9,0.9) to[out=90,in=0] ++(-0.9,0.9);

\draw (P2) ++(0.5,-0.3) node{$-a_1$};
\draw (P2) ++(3.1,-0.3) node{$-a_3$};
\draw (P2) ++(5.9,-0.3) node{$-a_{n-1}$};

\draw (P2) ++ (1.8,0.3) node{$\frac{a_2}{2}$};
\draw (P2) ++ (7.3,0.3) node{$\frac{a_n}{2}$};

\draw (P2) ++ (1.8,-0.9) node{$\frac{a_2}{2}$};
\draw (P2) ++ (7.3,-0.9) node{$\frac{a_n}{2}$};

\draw (P2) ++ (4.6,0.3) node{$\cdots$};
\draw (P2) ++ (4.6,-0.9) node{$\cdots$};

\draw[thick, densely dotted] (P2) ++(-0.7,-0.3) --++(0.7,0);
\draw[thick, densely dotted] (P2) ++(1,-0.3) --++(1.6,0);
\draw[thick, densely dotted] (P2) ++(3.6,-0.3) --++(1.7,0);
\draw[thick, densely dotted] (P2) ++(6.5,-0.3) --++(2.5,0);

\end{tikzpicture}

\caption{The diagram of a $2$-bridge knot $C(a_1 , \cdots, a_{n})$}
\label{fig:conway}
\end{figure}
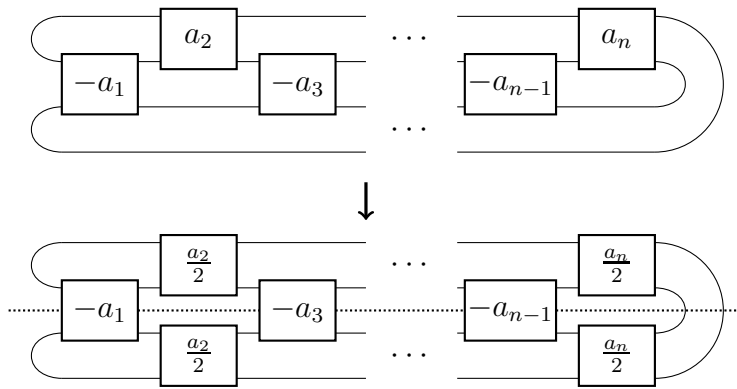

\end{document}